\documentclass[preprint,11pt, nopreprintline]{elsarticle}

\usepackage{amsmath}
\usepackage{amsfonts}
\usepackage{amssymb}
\usepackage{amsthm}

\usepackage{lineno}
\usepackage{multirow}
\usepackage{color}
\usepackage{subfigure}
\usepackage{algorithm}
\usepackage{algpseudocode}
\usepackage[margin=1in]{geometry}
\usepackage{tikz}
\usepackage[permil]{overpic}
\usepackage{hyperref}
\usepackage{cleveref}

\newtheorem{thm}{Theorem}[section]

\theoremstyle{definition}

\theoremstyle{remark}
\newtheorem{rem}[thm]{Remark}

\let\oldequation\equation
\let\oldendequation\endequation
\renewenvironment{equation}
  {\linenomathNonumbers\oldequation}
  {\oldendequation\endlinenomath}
  
\let\oldalign\align
\let\oldendalign\endalign
\renewenvironment{align}
  {\linenomathNonumbers\oldalign}
  {\oldendalign\endlinenomath}

\journal{Journal of Computational Physics}

\begin{document}

\begin{frontmatter}

\title{A correction function-based kernel-free boundary integral method for elliptic PDEs with implicitly defined interfaces}

\author[mysecondaryaddress]{Han Zhou}

\author[mysecondaryaddress]{Wenjun Ying\corref{mycorrespondingauthor}}
\cortext[mycorrespondingauthor]{Corresponding author}
\ead{wying@sjtu.edu.cn}

\address[mysecondaryaddress]{School of Mathematical Sciences, MOE-LSC and Institute of Natural Sciences, Shanghai Jiao Tong University, Minhang, Shanghai 200240, P. R. China}
\begin{abstract}

This work addresses a novel version of the kernel-free boundary integral (KFBI) method for solving elliptic PDEs with implicitly defined irregular boundaries and interfaces.
We focus on boundary value problems and interface problems, which are reformulated into boundary integral equations and solved with the matrix-free GMRES method.
In the KFBI method, evaluating boundary and volume integrals only requires solving equivalent but much simpler interface problems in a bounding box, for which fast solvers such as FFTs and geometric multigrid methods are applicable.
For the simple interface problem, a correction function is introduced for both the evaluation of right-hand side correction terms and the interpolation of a non-smooth potential function.
A mesh-free collocation method is proposed to compute the correction function near the interface.
The new method avoids complicated derivations of derivative jumps of the solution and is easy to implement, especially for the fourth-order method in three space dimensions.
Various numerical examples are presented, including challenging cases such as high-contrast coefficients, arbitrarily close interfaces and heterogeneous interface problems.
The reported numerical results verify that the proposed method is both accurate and efficient.

\end{abstract}


\end{frontmatter}


\section{Introduction}
Boundary value problems and interface problems of elliptic partial differential equations (PDEs) attract much attention due to their wide scientific and industrial applications, such as 
viscous incompressible flow \cite{quartapelle1993numerical, chorin1968numerical, CALHOUN2002231,RUSSELL2003177}, heat transfer \cite{bergman2011introduction, patankar2018numerical}, biomolecular electrostatics \cite{Davis1990, Honig1995},
electromagnetics \cite{rothwell2018electromagnetics, bondeson2012computational}, and many others.
In practical situations, domain boundaries and material interfaces are complex and even move with time, making it challenging to design accurate and efficient numerical methods for these problems.

Body-fitted discretization approaches, such as finite element methods \cite{Babuska1970, Guyomarch2009, Bramble1996, Wang2004}, approximate the computational domain with an unstructured mesh, which conforms to the geometry of boundaries and interfaces to achieve high-order accuracy.
However, it is always difficult and time-consuming to generate high-quality body-fitted meshes for complex geometries, especially when the boundary or interface moves substantially over time.
In addition, the linear systems generated from the discretization of the PDE on body-fitted meshes are less structured than those from a Cartesian grid, and fast solvers such as FFTs and geometric multigrid methods cannot be applied.

Immersed methods have been prevalent in recent decades, in which the complex boundary or interface is immersed into a fixed grid.
The pioneering work of immersed methods is the immersed boundary method (IBM) \cite{PESKIN1972252, Peskin2002, Peskin1977} that was initially proposed by C. S. Peskin for simulations of cardiac mechanics and blood flows.
In IBM, Peskin uses Lagrangian marker points on the boundary and regularized Dirac delta functions to approximate the singular force and spread it into the Eulerian grid.
The IBM is quite robust but is restricted to first-order accuracy due to the non-smoothness of the solution in the vicinity of the boundary.
Motivated by IBM, a number of immersed-type approaches have also been developed to improve the performance of conventional IBM.
Among them are the immersed interface method (IIM) \cite{Tan2008,Deng2003,Li1998,Leveque1994,Leveque1997,Li2001}, the ghost-fluid method (GFM) \cite{Fedkiw1999a,Fedkiw1999,Nguyen2001,Liu2000,Luo2006}, the matched interface and boundary (MIB) method \cite{ZHOU20061,Yu2007,Wang2015,Feng2019}, the correction function method (CFM) \cite{Marques2011,Marques2019,Marques2017}, and the Immersed Boundary Smooth Extension (IBSE) method \cite{Stein2017,Stein2016}.
The methods mentioned above are mainly based on finite difference discretizations.
Since the finite element method may provide more rigorous convergence analysis, similar ideas have also been used to develop finite element-based immersed methods, such as the extended finite element method (XFEM) \cite{Moes1999} and the immersed finite element method (IFEM) \cite{Li2003, Gong2008, Hou2005,Guo2020}.

The kernel-free boundary integral (KFBI) is a potential theory-based Cartesian grid method, which was initially proposed by W. Ying and C. S. Henriquez \cite{Ying2007} as an  extension of Mayo's method \cite{Mayo1984, Mayo1985, MCKENNEY1995348}. 
The KFBI method is also an immersed approach. Unlike traditional boundary integral methods (BIMs)/ boundary element methods (BEMs) \cite{Beale2004599, Beale2001,Greengard1998, Carrier1988,Carrier1988b, GREENGARD1987325,Liu1981, Ying2006,Klaseboer2004, TLUPOVA2009158}, 
layer and volume potentials are computed by solving equivalent but much simpler interface problems on a Cartesian grid, and the linear system can be efficiently solved with FFTs or geometric multigrid methods.
Therefore, the KFBI method has several attractive advantages: (a) no analytical expression of Green's function is needed for solving the boundary integral equation; (b) singular and nearly singular integrals are avoided; and (c) it can be applied to variable coefficient problems.
In the KFBI method, solving the constant coefficient interface problem is a fundamental building block.
In previous works \cite{Xie2020,Ying2014,Ying2013,Ying2007}, the simple interface problem is discretized with standard finite difference methods with a modified right-hand side.
The correction terms for the right-hand side are linear combinations of derivative jumps $[u], [u_x], [u_y], [u_z], [u_{xx}],\cdots$, which are computed by repeatedly taking tangential derivatives of the jump values and applying the local coordinate transformation.
The coordinate-transformation method for derivative jumps is accurate yet complicated when many derivative terms are needed, such as for high-order schemes and in three space dimensions \cite{Xie2020}.

In this work, we present a novel KFBI method that is both simple and accurate for two- and three-dimensional BVPs and interface problems.
Motivated by the correction function method (CFM)  \cite{Marques2011,Marques2019,Marques2017}, we introduce a correction function in the vicinity of the interface to derive correction terms of the right-hand side for the constant coefficient interface problem.
In order to solve the local Cauchy problem for the correction function, we propose a mesh-free collocation method based on an overlapping surface decomposition for the interface.
Unlike the original CFM \cite{Marques2011,Marques2019,Marques2017}, no surface quadrature is required since the collocation method works with the strong form of the Cauchy problem.
The overlapping surface decomposition representation of the interface also provides a good choice of collocation points such that the resulting collocation problem is accurate and stable.
Another property of the collocation method is that the discrete system of the collocation problem is a square one and can be solved accurately, which is different from the original CFM in that the linear system is overdetermined and needs to be solved in the least-square sense.
The new approach for the constant coefficient interface problem is built into the KFBI framework to accommodate elliptic BVPs and more general interface problems.
The resulting method is named the correction function-based KFBI method.

The paper is organized as follows.
The governing equations and their boundary integral equations are described in section \ref{sec:problem} and \ref{sec:bie}.
In section \ref{sec:kfbi}, the main idea of the KFBI method is described.
The details of the numerical method for the constant coefficient interface problem are described in section \ref{sec:crc}.
The algorithm is summarized in section \ref{sec:algo}.
In section \ref{sec:res}, numerical results demonstrating the method with examples are presented.
Finally, we discuss the improvement and advantages of the proposed method in section \ref{sec:dis}.

\section{Governing equations}\label{sec:problem}

\subsection{Boundary value problem}
Let $\Omega \subset \mathbb{R}^d, d = 2,3$ be a complex domain with smooth boundary $\Gamma = \partial\Omega$, as illustrated in \Cref{fig:domain-1}. 
The BVP of an elliptic PDE is given by
\begin{equation}\label{eqn:BVP}
    \nabla \cdot (\sigma \nabla u)- \kappa u = f,\quad \text{in } \Omega,
\end{equation}
subject to either the Dirichlet boundary condition or the Neumann boundary condition
\begin{equation}\label{eqn:BVP-bc}
    u = g_D,\quad \text{or}\quad \sigma\partial_{\boldsymbol{n}}u = g_N,\quad \text{on } \Gamma,
\end{equation}
where $\sigma >0$ is the diffusivity and $\kappa\geq 0$ is the reaction coefficient.
In this paper, we assume that $\sigma$ and $\kappa$ are constants.

\subsection{Interface problem}
\begin{figure}[htbp]
    \centering
    \subfigure[]{\includegraphics[width=0.3\textwidth]{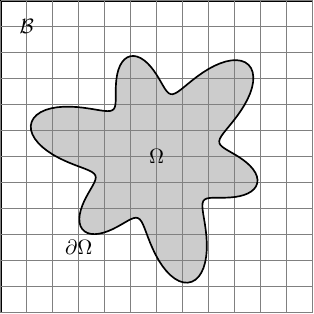}\hspace{.2in}\label{fig:domain-1}}
    \subfigure[]{\includegraphics[width=0.3\textwidth]{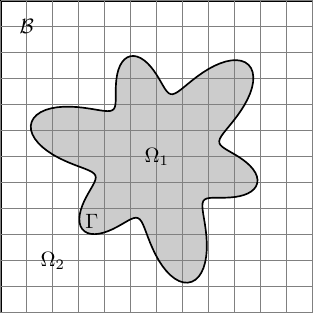}\hspace{.2in}\label{fig:domain-2}}
    \caption{A schematic of the (a) boundary value problem and (b) interface problem. Irregular domains and interfaces are embedded into a larger bounding box, in which a uniform Cartesian grid is used for computation.}
    \label{fig:domains}
\end{figure}

Let $\Gamma\subset \mathbb{R}^d, d = 2,3$ be a sharp interface that separates a larger domain $\mathcal{B}\subset \mathbb{R}^d$ into two subdomains $\Omega_1$ and $\Omega_2$, as illustrated in \Cref{fig:domain-2}.
The interface problem of an elliptic PDE is given by
\begin{equation}\label{eqn:IFP}
    \nabla \cdot (\sigma_i \nabla u) - \kappa_i u = f_i, \quad \text{in } \Omega_i,\quad i = 1, 2,
\end{equation}
subject to two interface jump conditions
\begin{equation}\label{eqn:IFP-bc}
    [u] = g_1,\quad [\sigma \partial_{\boldsymbol{n}}u] = g_2, \quad \text{on } \Gamma,
\end{equation}
and a homogeneous Dirichlet boundary condition on the outer boundary
\begin{equation}
\label{eqn:diri-bc}
    u = 0,\quad\text{on }\partial\mathcal{B},
\end{equation}
where $\sigma_i >0,i=1,2$ are diffusivities and $\kappa_i\geq 0,i=1,2$ are reaction coefficients.
Similarly, we only consider the case in which $\sigma_i$ and $\kappa_i$ are constants.
Note that the Dirichlet boundary condition \eqref{eqn:diri-bc} is chosen only for simplicity, since the treatment for boundary conditions on $\partial\mathcal{B}$ only depends on the finite difference scheme and is much simpler.
Different boundary conditions, such as Neumann and periodic ones, can also be used.

Here, the boundary/interface $\Gamma$ is assumed to be implicitly defined as the zero level set of a function.
In the case that $\Gamma$ is defined by a parametric surface or spline, it can also be transformed into an implicit form.

\section{Boundary integral equations}\label{sec:bie}
Both the boundary value problem \eqref{eqn:BVP}--\eqref{eqn:BVP-bc} and the interface problem \eqref{eqn:IFP}--\eqref{eqn:IFP-bc} are solved by reformulating them as boundary integral equations.

\subsection{Boundary value problem}
Let $G(\boldsymbol{q}, \boldsymbol{p})$ be Green's function such that for each fixed $\boldsymbol{p}\in\mathcal{B}$, 
\begin{equation}
    \begin{aligned}
    \nabla_{\boldsymbol{q}}\cdot (\sigma(\boldsymbol{q}) \nabla_{\boldsymbol{q}} G(\boldsymbol{q},\boldsymbol{p})) - \kappa(\boldsymbol{q})G(\boldsymbol{q};\boldsymbol{p}) &= \delta(\boldsymbol{q} - \boldsymbol{p}), \quad &\text{in } \mathcal{B},\\
    G(\boldsymbol{q},\boldsymbol{p})& = 0, \quad &\text{on }\partial \mathcal{B}.
    \end{aligned}
\end{equation}
Let $\varphi, \psi$ be two density functions.
Define the single layer, double layer, adjoint double layer and hyper-singular integrals, respectively, by
\begin{align}
    \mathcal{S}\psi(\boldsymbol{p}) &= \int_{\Gamma}G(\boldsymbol{q};\boldsymbol{p})\psi(\boldsymbol{q}) \,
    d\boldsymbol{s}_{\boldsymbol{q}},\quad& \boldsymbol{p}\in \Gamma, \\
    \mathcal{K}\varphi(\boldsymbol{p}) &= \int_{\Gamma}\sigma(\boldsymbol{q})  \dfrac{\partial G(\boldsymbol{q};\boldsymbol{p})}{\partial \boldsymbol{n}_{\boldsymbol{q}}}\varphi(\boldsymbol{q}) \,
    d\boldsymbol{s}_{\boldsymbol{q}},\quad &\boldsymbol{p}\in \Gamma, \\
    \mathcal{K}^{\prime}\psi(\boldsymbol{p}) &= \int_{\Gamma} \sigma(\boldsymbol{p})\dfrac{\partial G(\boldsymbol{q};\boldsymbol{p})}{\partial \boldsymbol{n}_{\boldsymbol{p}}}\psi(\boldsymbol{q}) \,
    d\boldsymbol{s}_{\boldsymbol{q}},\quad &\boldsymbol{p}\in \Gamma, \\
    \mathcal{D}\varphi(\boldsymbol{p}) &= \int_{\Gamma}\sigma(\boldsymbol{p})\sigma(\boldsymbol{q}) \dfrac{\partial^2 G(\boldsymbol{q};\boldsymbol{p})}{\partial \boldsymbol{n}_{\boldsymbol{q}}\partial \boldsymbol{n}_{\boldsymbol{p}}}\varphi(\boldsymbol{q}) \,
    d\boldsymbol{s}_{\boldsymbol{q}},\quad& \boldsymbol{p}\in \Gamma.
\end{align}
Define the volume integrals by
\begin{align}
    \mathcal{G}f(\boldsymbol{p}) &= \int_{\Omega} G(\boldsymbol{q};\boldsymbol{p})f(\boldsymbol{q})\, d\boldsymbol{q},\quad \boldsymbol{p}\in \Gamma\label{eqn:vol-int}, \\
    \partial_{\mathbf{n}} \mathcal{G}f(\boldsymbol{p}) &= \dfrac{\partial}{\partial \mathbf{n}_{\mathbf{p}}} \int_{\Omega} G(\boldsymbol{q};\boldsymbol{p})f(\boldsymbol{q})\, d\boldsymbol{q},\quad \boldsymbol{p}\in \Gamma\label{eqn:vol-int-norm}.
\end{align}

The Dirichlet and Neumann BVPs \eqref{eqn:BVP}--\eqref{eqn:BVP-bc} can be reformulated, respectively, as the boundary integral equations
\begin{align}
    (\dfrac{1}{2} + \mathcal{K})\varphi &= g_D -\mathcal{G}f,\label{eqn:bie-D}\\
    (\dfrac{1}{2} -\mathcal{K}^{\prime}) \psi &= g_N- \sigma\partial_{ \boldsymbol{n}}\mathcal{G}f,\label{eqn:bie-N}
\end{align}
which are both Fredholm integral equations of the second kind and well-conditioned.

\subsection{Interface problem}
Let $G_i(\boldsymbol{q}, \boldsymbol{p}), i=1,2$ be Green's functions such that for each fixed $\boldsymbol{p}\in\mathcal{B}$, 
\begin{equation}
    \begin{aligned}
    \nabla_{\boldsymbol{q}}\cdot (\sigma_i(\boldsymbol{q}) \nabla_{\boldsymbol{q}} G_i(\boldsymbol{q},\boldsymbol{p})) - \kappa_i(\boldsymbol{q})G_i(\boldsymbol{q};\boldsymbol{p}) &= \delta(\boldsymbol{q} - \boldsymbol{p}), \quad &\text{in } \mathcal{B},\\
    G_i(\boldsymbol{q},\boldsymbol{p})& = 0, \quad &\text{on }\partial \mathcal{B}.
    \end{aligned}
\end{equation}
Similarly, we can define the single layer, double layer, adjoint double layer, hyper-singular and volume integral operators $\mathcal{S}_i, \mathcal{K}_i, \mathcal{K}^{\prime}_i, \mathcal{D}_i, \mathcal{G}_i, \partial_{\mathbf{n}}\mathcal{G}_i, i = 1,2$.
By introducing two unknown density functions $\varphi = u_1$ and $\psi = \sigma_2\partial_{\boldsymbol{n}}u_2$, the interface problem \eqref{eqn:IFP} can be reformulated as a system of boundary integral equations
\begin{equation}\label{eqn:bie-ifp-1}
\begin{aligned}
    \varphi - (\mathcal{K}_1 - \mathcal{K}_2)\varphi + (\mathcal{S}_1 - \mathcal{S}_2)\psi &= \dfrac{1}{2} g_1 + \mathcal{G}_1 f_1 + \mathcal{G}_2 f_2 + \mathcal{K}_2g_1 - \mathcal{S}_1 g_2,\\  \psi  - (\mathcal{D}_1 - \mathcal{D}_2)\varphi + (\mathcal{K}^{\prime}_1 - \mathcal{K}^{\prime}_2)\psi &= -\dfrac{1}{2} g_2 + \sigma_1\partial_{\boldsymbol{n}}\mathcal{G}_1 f_1 + \sigma_2\partial_{\boldsymbol{n}}\mathcal{G}_2 f_2 + \mathcal{D}_2 g_1 - \mathcal{K}^{\prime}_1 g_2.
\end{aligned}
\end{equation}

In the case of $\kappa_i = 0, i=1,2$ or $\sigma_1/\sigma_2 = \kappa_1/\kappa_2$, dividing the two equations in \eqref{eqn:IFP} by $\sigma_i$, respectively, yields
\begin{equation}\label{eqn:bie-ifp-2}
    \Delta u - \Tilde{\kappa} u =
    \left \{
    \begin{aligned}
    f_1 / \sigma_1, \quad \text{in }\Omega_1, \\
    f_2 / \sigma_2, \quad \text{in }\Omega_2,
    \end{aligned}
    \right .
\end{equation}
where $\Tilde{\kappa}=0$ or $\Tilde{\kappa}=\kappa_1/\sigma_1=\kappa_2/\sigma_2$.
With $\psi = [\partial_{\boldsymbol{n}}u]$, we may also obtain a simpler boundary integral equation
\begin{equation}
    \dfrac{1}{2}\psi + \mu \mathcal{K}^{\prime}\psi = \dfrac{ g_2}{\sigma_1+\sigma_2} + \mu (\mathcal{D}g_1 + \partial_{\boldsymbol{n}}\mathcal{G}f),
\end{equation}
where $\mu = (\sigma_2 - \sigma_1)/(\sigma_2 + \sigma_1)\in(-1,1)$ and $\mathcal{K}, \mathcal{D}$ and $\mathcal{G}$ are the integral operators associated with the Green's function of the operator $\Delta-\Tilde{\kappa}$.
One may refer to \cite{Ying2014} for detailed derivations of the boundary integral equations.

\section{Kernel-free boundary integral method}\label{sec:kfbi}
In the kernel-free boundary integral method, values of the boundary and volume integrals at the boundary or interface $\Gamma$ are not evaluated with quadrature methods.
Instead, they are evaluated by solving equivalent but much simpler interface problems for boundary and volume potentials.

For Green's function $G(\boldsymbol{q},\boldsymbol{p})$, which is associated with the elliptic operator $\sigma\Delta - \kappa$, define the single layer potential $-S\psi$, the double layer potential $D\varphi$ and the Newtonian potential $Nf$ by
\begin{equation}
    \begin{aligned}
    -S\psi(\boldsymbol{p}) &=- \int_{\Gamma}G(\boldsymbol{q},\boldsymbol{p})\psi(\boldsymbol{q}) \,
    d\boldsymbol{s}_{\boldsymbol{q}},\quad &\boldsymbol{p}\in\mathcal{B}, \\
    D\varphi(\boldsymbol{p}) &= \int_{\Gamma}\sigma(\boldsymbol{q})  \dfrac{\partial G(\boldsymbol{q},\boldsymbol{p})}{\partial \boldsymbol{n}_{\boldsymbol{q}}}\varphi(\boldsymbol{q}) \,
    d\boldsymbol{s}_{\boldsymbol{q}},\quad &\boldsymbol{p}\in\mathcal{B},\\
    Nf(\boldsymbol{p}) &= \int_{\Omega} G(\boldsymbol{q};\boldsymbol{p})f(\boldsymbol{q})\, d\boldsymbol{q},\quad& \boldsymbol{p}\in \mathcal{B}.
    \end{aligned}
\end{equation}
Then the boundary integrals $\mathcal{S}\psi$, $\mathcal{K}\varphi$, $\mathcal{K}^{\prime}\psi$, $\mathcal{D}\varphi$ and the volume integrals $\mathcal{G}f, \partial_{\mathbf{n}}\mathcal{G}f$ coincide with boundary values or normal derivatives of the potentials $S\psi$, $D\varphi$ and $Nf$.
The above three potential functions are not smooth at $\Gamma$ and, by classical potential theory, satisfy equivalent interface problems (see \cite{Ying2007,Xie2020}).
The equivalent interface problems for the single layer potential $-S\psi$, the double layer potential $D\varphi$ and the Newton potential $Nf$ can be unified as
\begin{equation}\label{eqn:equi-ifp}
    \left \{
    \begin{aligned}
    &\nabla\cdot(\sigma\nabla u) - \kappa u = F, \quad & \text{in }\mathcal{B}\setminus \Gamma, \\
    &[u] = \Phi, \quad &\text{on }\Gamma, \\
    &[\partial_{\boldsymbol{n}} u] = \Psi, \quad &\text{on }\Gamma, \\
    &u = 0,\quad &\text{on }\partial\mathcal{B}.
    \end{aligned}
    \right .
\end{equation}
The functions $\Phi$, $\Psi$ and $F$ are specified for each potential by
\begin{itemize}
    \item $-S\psi$: $\Phi = F = 0$, $\Psi = \psi$.
    \item $D\varphi$: $\Phi = \varphi$, $\Psi= F=0$. 
    \item $Nf$:  $\Phi = \Psi = 0$. $F$ is an arbitrary extension of $f$ to the whole box $\mathcal{B}$. For simplicity, we set the extended value as zero.
\end{itemize}

Once the interface problem \eqref{eqn:equi-ifp} is solved for the potentials, the boundary integrals $\mathcal{S}\psi$, $\mathcal{K}\varphi$, $\mathcal{K}^{\prime}\psi$, $\mathcal{D}\varphi$ and the volume integrals $\mathcal{G}f, \partial_{\mathbf{n}}\mathcal{G}f$ can be obtained from the grid data of these potentials with an interpolation method.

\section{Equivalent simple interface problem}\label{sec:crc}
Solving the constant coefficient interface problem \eqref{eqn:equi-ifp} is an essential part of the KFBI method.
For simplicity, we drop the constant $\sigma$ and proceed with the following problem
\begin{equation}\label{eqn:con-ifp}
    \left \{
    \begin{aligned}
    &\Delta u- \kappa u = f, \quad &\text{in } \mathcal{B}\setminus \Gamma,\\
    &[u] = a, \quad &\text{on }\Gamma, \\
    &[\partial_{\boldsymbol{n}} u] = b, \quad &\text{on }\Gamma, \\
    &u = 0,\quad &\text{on }\partial\mathcal{B}.
    \end{aligned}
    \right .
\end{equation}
where $a$, $b$ and $f$ are given data and $\kappa$ is a constant.
The right-hand side $f$ is possibly discontinuous across the interface $\Gamma$.
The constant coefficient interface problem is a much simpler case of the more general interface problem \eqref{eqn:IFP}.

\subsection{Interface representation}

In this work, the interface $\Gamma$ is implicitly defined by the level set function $H$ in the following way
\begin{equation}
    \Gamma = \{\mathbf{x} \in \mathbb{R}^d | H(\mathbf{x}) = 0\}, \quad \text{for }d = 2,3.
\end{equation}
We assume the level set function $H$ is at least $C^4$ (for the fourth-order method) and $|\nabla H|>c_0$ for some $c_0 > 0$ near the interface $\Gamma$.
The level set function allows us to easily determine the intersection points of the surface with grid lines. For instance, if we have an intersection point on the line segment between two grid nodes $\mathbf{x}_1$ and $\mathbf{x}_2$, we expect the values $H(\mathbf{x}_1)$ and $H(\mathbf{x}_2)$ to have opposite signs. By solving the scalar algebraic equation for $t$ as follows:
\begin{equation}
    H(t\mathbf{x}_1 + (1-t)\mathbf{x}_2) = 0, \quad t \in [0,1],
\end{equation}
using methods such as Newton's method or the bisection method, one can obtain the coordinates of the intersection point.
To compute the unit outward normal at a surface point, we utilize the gradient of the level set function. The unit outward normal vector $\mathbf{n}(\mathbf{x})$ at a point $\mathbf{x}$ on the surface is given by:
\begin{equation}
    \mathbf{n}(\mathbf{x}) = \dfrac{\nabla H(\mathbf{x})}{|\nabla H(\mathbf{x})|}.
\end{equation}

\begin{figure}[htbp]
    \centering
    \subfigure[]{\includegraphics[width=0.3\textwidth]{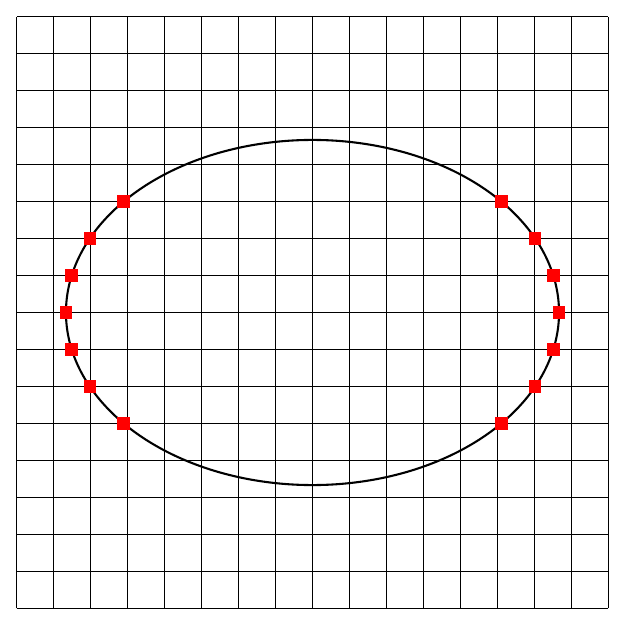}}
    \subfigure[]{\includegraphics[width=0.3\textwidth]{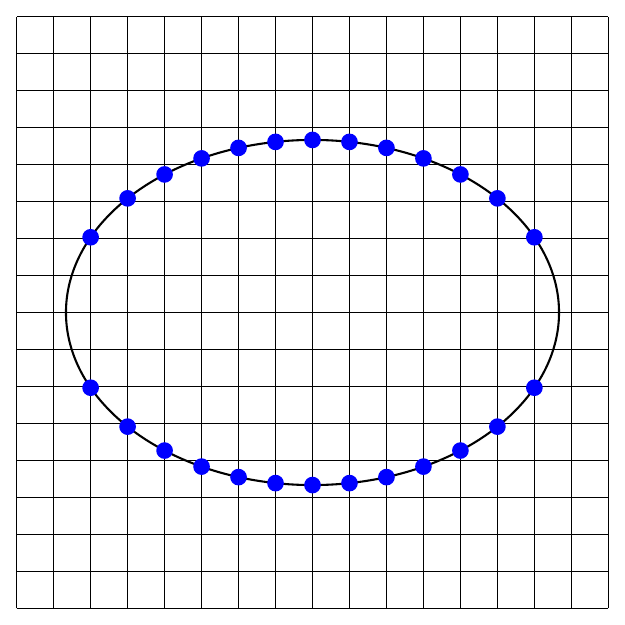}}
    \subfigure[]{\includegraphics[width=0.3\textwidth]{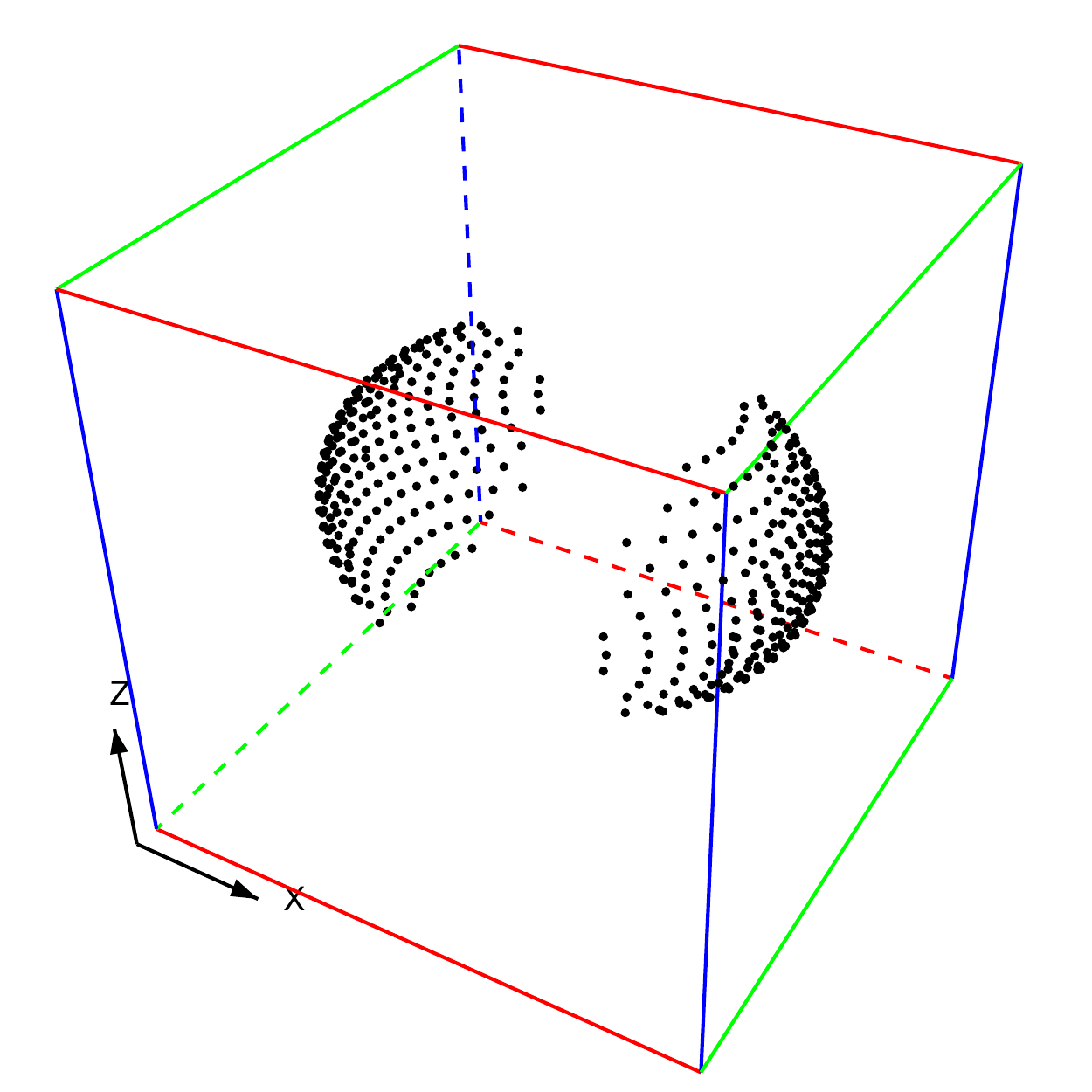}}
    \caption{Illustrations of surface points: (a) points in $\Gamma_1^h$ (in 2D); (b) points in $\Gamma_2^h$ (in 2D); (c) points in $\Gamma_1^h$ (in 3D).}
    \label{fig:surface-points}
\end{figure}

The method described in \cite{Ying2013} is employed in this work for the representation of the interface using only a subset of intersection points. For each $r = 1,\ldots,d$, $\mathbf{e}_r$ represents the $r$-th Cartesian basis vector in $\mathbb{R}^d$, and $\alpha \in (\cos^{-1}(1/\sqrt{d}), \pi/2)$ is a fixed angle. We define the subset:
\begin{equation}
    \Gamma_r = \{\mathbf{x} \in \Gamma : |\mathbf{n}(\mathbf{x}) \cdot \mathbf{e}_r | > \cos\alpha\},
\end{equation}
which forms an overlapping surface decomposition of $\Gamma$. The discrete representation of the interface $\Gamma_r$ only considers the intersection points between $\Gamma_r$ and the grid lines aligned with the $\mathbf{e}_r$ direction (refer to \Cref{fig:surface-points}). We denote the set of these intersection points as $\Gamma_r^h$. The union of all sets $\Gamma_r^h$ for $r = 1,\ldots,d$ is denoted as $\Gamma^h$, which represents the discrete set of points used to approximate $\Gamma$ and allocate surface degrees of freedom.
For more detailed information on the surface discretization algorithm, please refer to \cite{Ying2013}.

With the help of the overlapping surface decomposition-based discretization, the interface $\Gamma$ can be locally parameterized by a reference coordinate plane.
Candidate reference planes are 
\begin{align}
    \Pi_i:\{(x_1, x_2, \cdots,x_d)\in\mathbb{R}^d|x_i = 0\}, \quad i=1,\cdots,d,
\end{align}
for $d = 2 \text{ or } 3$.
Suppose that, at a point $\boldsymbol{x}\in\Gamma$, the $i$-th component of the local normal $\boldsymbol{n}(\boldsymbol{x})$ has the largest absolute value. Then we choose $\Pi_i$ as the reference plane of $\Gamma$ near $\boldsymbol{x}$.
In such a way, the interpolation stencils on $\Gamma$ can be easily found with the help of the Cartesian grid on the reference plane.
Numerical integration and interpolation on $\Gamma$ can be done in a way similar to those on a planar domain.
We remark that, in principle, the Cartesian grid used for the representation of $\Gamma$ is not necessarily the same as the one used for solving PDEs.
In this work, we use the same Cartesian grid only for simplicity.

\subsection{Corrected finite difference scheme}

For simplicity, the bounding box is assumed to be a unit cube, i.e., $\mathcal{B} = (0,1)^3$.
Given a positive integer $N$, the domain $\mathcal{B}$ is uniformly partitioned into a Cartesian grid with mesh parameter $h = 1/N$.
Let $P_{i, j, k}$ denote the grid node $(x_i, y_j, z_k), i, j, k = 0,1, \cdots, N$, where $x_i=ih$, $y_j = jh$ and $z_k = kh$ are node coordinates.
For an irregular domain $\Omega\subset\mathcal{B}$, the interior and exterior grid nodes are defined as $\Omega_h$ and $\Omega_h^C$, respectively,
\begin{equation}
\begin{aligned}
     \Omega_h &= \{P_{i,j,k} | (x_i,y_j,z_k) \in \Omega,\quad i,j,k = 1,\cdots,N-1 \},\\
    \Omega_h^C &= \{P_{i,j,k} | (x_i,y_j,z_k) \in \mathcal{B}\setminus\Omega,\quad i,j,k = 1,\cdots,N-1 \}.  
\end{aligned}
\end{equation}

In the absence of interfaces, it is known that the following two compact finite difference schemes \eqref{eqn:9p-cds} and \eqref{eqn:27p-cds} are fourth-order accurate for 2D and 3D cases, respectively.
\begin{equation}\label{eqn:9p-cds}
    -(\dfrac{10}{3h^2}+\dfrac{2}{3}\kappa)u_{i,j} 
    + (\dfrac{2}{3h^2} - \dfrac{1}{12}\kappa)\sum_{\substack{r,s\in\{-1,0,1\}\\ |r|+|s|=1}} u_{i+r,j+s} + \dfrac{1}{6h^2}\sum_{\substack{r,s\in\{-1,0,1\}\\ |r|+|s|=2}} u_{i+r,j+s}
     = f_{i,j}+\dfrac{h^2}{12}\Delta f_{i,j}.
\end{equation}
\begin{equation}\label{eqn:27p-cds}
\begin{aligned}
    &-(\dfrac{25}{6h^2}+\dfrac{1}{2}\kappa)u_{i,j,k} 
    + (\dfrac{5}{12h^2} - \dfrac{1}{12}\kappa)\sum_{\substack{r,s,t\in\{-1,0,1\}\\ |r|+|s|+|t|=1}} u_{i+r,j+s,k+t}
    + \dfrac{1}{8h^2}\sum_{\substack{r,s,t\in\{-1,0,1\}\\ |r|+|s|+|t|=2}} u_{i+r,j+s,k+t}\\
    &+\dfrac{1}{48h^2}\sum_{\substack{r,s,t\in\{-1,0,1\}\\ |r|+|s|+|t|=3}} u_{i+r,j+s,k+t}
     = f_{i,j,k}+\dfrac{h^2}{12}\Delta f_{i,j,k}.
\end{aligned}
\end{equation}
The two schemes are adopted to derive the corresponding corrected finite difference schemes for the interface problem \eqref{eqn:con-ifp}.
We write the finite difference schemes in the general form
\begin{equation}\label{eqn:FDS-gene}
    \sum_{P_{i+r,j+s,k+t}\in \mathcal{S}_{i,j,k}}c_{r,s,t}u_{i+r, j+s,k+t} = F_{i,j,k},
\end{equation}
where $c_{r,s,t}$ is the coefficient of $u_{i+r,j+s,k+t}$, $F_{i,j,k}$ is the right-hand side of the finite difference equation and $\mathcal{S}_{i,j,k}$ is the node set that contains all grid nodes with $c_{r,s,t} \neq 0$ at $P_{i,j,k}$.
Then we define regular nodes $\mathcal{R}_h$ and irregular nodes $\mathcal{I}_h$ as follows,
\begin{align}
    \mathcal{R}_h = \{ P_{i,j,k} | S_{i,j,k}\cap \Omega_h = \emptyset \text{ or } S_{i,j,k}\cap \Omega_h^C = \emptyset \}, \\
    \mathcal{I}_h = \{ P_{i,j,k} | S_{i,j,k}\cap \Omega_h \neq \emptyset \text{ and } S_{i,j,k}\cap \Omega_h^C \neq \emptyset \},
\end{align}

At irregular nodes, since the finite difference approximation is taken across the discontinuity at the interface, large local truncation errors may occur and result in inaccurate or even divergent results.
Precisely, let $\mathcal{A}_h$ denote the difference operator in the finite difference scheme \eqref{eqn:FDS-gene}.
Suppose the local truncation error is on the order of $\mathcal{O}(h^p)$ at a regular node.
The local truncation error at an irregular node is given by
\begin{equation}\label{eqn:LTE}
    \begin{aligned}
    &E_h(x_i, y_j, z_k) = \mathcal{A}_h u(x_i, y_j, z_k) - F_{i,j,k} \\=
    &\left \{
    \begin{aligned}
    &\sum_{P_{i+r,j+s,k+t}\in \Omega_h\cap\mathcal{S}_{i,j,k}}c_{r,s,t}u^+(x_{i+r},y_{j+s},z_{k+t})\\
    &+\sum_{P_{i+r,j+s,k+t}\in \Omega_h^C\cap\mathcal{S}_{i,j,k}}c_{r,s,t}u^-(x_{i+r},y_{j+s},z_{k+t})- F_{i,j,k},\quad P_{i,j,k}\in\Omega_h,\\
   &\sum_{P_{i+r,j+s,k+t}\in \Omega_h^C\cap\mathcal{S}_{i,j,k}}c_{r,s,t}u^-(x_{i+r},y_{j+s},z_{k+t}) \\
    &+\sum_{P_{i+r,j+s,k+t}\in \Omega_h\cap\mathcal{S}_{i,j,k}}c_{r,s,t}u^+(x_{i+r},y_{j+s},z_{k+t})- F_{i,j,k},\quad P_{i,j,k}\in\Omega^C_h,\\
    \end{aligned}
    \right . \\=
    &\left \{
    \begin{aligned}
    &\sum_{P_{i+r,j+s,k+t}\in \mathcal{S}_{i,j,k}}c_{r,s,t}u^+(x_{i+r},y_{j+s},z_{k+t}) \\
    &+\sum_{P_{i+r,j+s,k+t}\in \Omega_h^C\cap\mathcal{S}_{i,j,k}}c_{r,s,t}(u^- -u^+)(x_{i+r},y_{j+s},z_{k+t})- F_{i,j,k},\quad P_{i,j,k}\in\Omega_h,\\
    &\sum_{P_{i+r,j+s,k+t}\in \mathcal{S}_{i,j,k}}c_{r,s,t}u^-(x_{i+r},y_{j+s},z_{k+t})\\
    & +\sum_{P_{i+r,j+s,k+t}\in \Omega_h\cap\mathcal{S}_{i,j,k}}c_{r,s,t}(u^+ - u^-)(x_{i+r},y_{j+s},z_{k+t})- F_{i,j,k},\quad P_{i,j,k}\in\Omega^C_h,\\
    \end{aligned}
    \right . \\=
    &\left \{
    \begin{aligned}
    \dfrac{1}{h^2} \sum_{P_{i+r,j+s,k+t}\in\Omega_h^C\cap\mathcal{S}_{i,j,k}}c_{r,s,t}(u^- - u^+)(x_{i+r},y_{j+s},z_{k+t}) +  \mathcal{O}(h^p),\quad P_{i,j,k}\in\Omega_h,\\
     \dfrac{1}{h^2} \sum_{P_{i+r,j+s,k+t}\in\Omega_h\cap\mathcal{S}_{i,j,k}}c_{r,s,t}(u^+ - u^-)(x_{i+r},y_{j+s},z_{k+t}) + \mathcal{O}(h^p),\quad P_{i,j,k}\in\Omega^C_h,\\
    \end{aligned}
    \right .
    \end{aligned}
\end{equation}
where $u^+$ and $u^-$ are two smooth functions that coincide with $u$ in the domain $\Omega$ and $\Omega^C$, respectively.
It can be found that the leading term in the local truncation error at an irregular node is on the order of $\mathcal{O}(h^{-2})$, which is not acceptable for the sake of accuracy.
The problem can be fixed by including the leading terms of the local truncation error, as correction terms, into the final finite difference equations.
Define the correction function $C(\boldsymbol{x}) = u^+(\boldsymbol{x}) - u^-(\boldsymbol{x})$.
Then, the corrected finite difference scheme can be written as
\begin{equation}\label{eqn:corrected-FDS}
    \sum_{P_{i+r,j+s,k+t}\in \mathcal{S}_{i,j,k}}c_{r,s,t}u_{i+r, j+s,k+t} = F_{i,j,k}+C_{i,j,k}.
\end{equation}
where the correction term $C_{i,j,k}$ is given by
\begin{equation}\label{eqn:crc-formula}
    C_{i,j,k} = \left \{
    \begin{aligned}
    &0,\quad &P_{i,j,k}\in\mathcal{R}_h, \\
    &-\dfrac{1}{h^2} \sum_{P_{i+r,j+s,k+t}\in\Omega_h^C\cap\mathcal{S}_{i,j,k}}c_{r,s,t}C(x_{i+r},y_{j+s},z_{k+t}),\quad &P_{i,j,k}\in\Omega_h\cap\mathcal{I}_h,\\
    &\dfrac{1}{h^2} \sum_{P_{i+r,j+s,k+t}\in\Omega_h\cap\mathcal{S}_{i,j,k}}c_{r,s,t}C(x_{i+r},y_{j+s},z_{k+t}),\quad &P_{i,j,k}\in\Omega^C_h\cap\mathcal{I}_h.\\
    \end{aligned}
    \right .
\end{equation}

\begin{rem}
If exact values of the correction function $C(\boldsymbol{x})$ are given, then the local truncation error of the corrected finite difference scheme \eqref{eqn:corrected-FDS} is on the order of $\mathcal{O}(h^p)$ at each node. 
However, it happens only when the interface coincides with grid nodes and $C(\boldsymbol{x})$ equals the Dirichlet jump condition $[u]$.
In practice, approximate values of the correction function $C(\boldsymbol{x})$ are used.
For the fourth-order method in this work, the correction function only needs to be approximated with an error on the order of $\mathcal{O}(h^5)$ such that the local truncation error becomes $\mathcal{O}(h^3)$ at irregular nodes and $\mathcal{O}(h^4)$ elsewhere.

\end{rem}

\begin{figure}[htbp]
\centering
\begin{overpic}[scale=0.7]{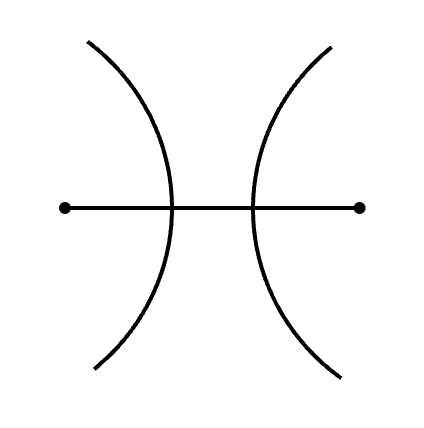}
\put(130,450){$P_{i,j}$}
\put(800,450){$P_{i+r,j+s}$}
\put(200,700){$\Gamma^{(1)}$}
\put(700,700){$\Gamma^{(2)}$}
\put(200,350){$\Omega^{(1)}$}
\put(700,350){$\Omega^{(2)}$}
\put(450,250){$\Omega^{(0)}$}
\end{overpic}
\caption{An illustration of a line segment intersecting two interfaces.}
\label{fig:double-cross}
\end{figure}
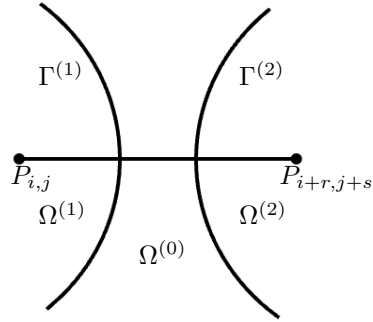

\begin{rem}
If two interfaces are arbitrarily close, the line segment between two grid nodes may intersect interfaces more than once (see \Cref{fig:double-cross}).
Let $u^{(i)}$ be the restrictions of the piecewise smooth solution $u$ to $\Omega^{(i)}$ for $i=0,1,2$.
Denote by $C^{(1)} = u^{(1)}-u^{(0)}$ and $C^{(2)}=u^{(2)}-u^{(0)}$ two correction functions that are computed near $\Gamma^{(1)}$ and $\Gamma^{(2)}$.
For the correction term $C_{i,j}$, the value $C(x_{i+r},y_{j+s})$ in \eqref{eqn:crc-formula} is computed by
\begin{equation}
\begin{aligned}
    C(x_{i+r},y_{j+s}) &= u^{(1)}(x_{i+r},y_{j+s})-u^{(2)}(x_{i+r},y_{j+s}) \\
    &= u^{(1)}(x_{i+r},y_{j+s}) - u^{(0)}(x_{i+r},y_{j+s}) + u^{(0)}(x_{i+r},y_{j+s}) - u^{(2)}(x_{i+r},y_{j+s}) \\
    &= C^{(1)}(x_{i+r},y_{j+s}) - C^{(2)}(x_{i+r},y_{j+s}).        
\end{aligned}
\end{equation}
This is simply adding and subtracting a middle term and is similar to the technique used in \cite{Ying2018}.
\end{rem}

\subsection{Local Cauchy problem}

Suppose $\Gamma$ is sufficiently smooth and the right-hand side $f$ is also piecewise smooth.
Denote by $\Omega_{\Gamma}$ a narrow band around $\Gamma$ that covers all irregular nodes.
Let $f^+$ and $f^-$ be two smooth extension functions of $f$ in $\Omega_{\Gamma}$ from two different sides $\Omega_i$ and $\Omega_e$, respectively.
Then the function $\Tilde{f} = f^+ - f^-$ is also smooth in $\Omega_{\Gamma}$.
The smoothness of $\Tilde{f}$ is relevant to the accuracy of $C(x)$, see \Cref{rem:f-smooth}.
Notice that the correction function $C(\boldsymbol{x})$ satisfies the Cauchy problem 
\begin{equation}\label{eqn:cau-prb}
    \begin{aligned}
    \Delta C(\boldsymbol{x}) - \kappa C(\boldsymbol{x})&= \Tilde{f}(\boldsymbol{x}),\quad &\boldsymbol{x}\in\Omega_{\Gamma},\\
    C(\boldsymbol{x}) &= a(\boldsymbol{x}),\quad &\boldsymbol{x}\in\Gamma,\\
    \partial_{\boldsymbol{n}} C(\boldsymbol{x}) &= b(\boldsymbol{x}), \quad &\boldsymbol{x}\in \Gamma.
    \end{aligned}
\end{equation}
The Cauchy problem is known to be ill-posed in the sense of Hadamard: small perturbations in the boundary data grow exponentially away from the boundary, making it difficult to obtain a global numerical solution.
Since the correction function is only required at irregular nodes that are close to the boundary $\Gamma$, we are only interested in the local solution of the Cauchy problem.
In that case, numerical errors can be bounded from above.
The localness of the Cauchy problem also suggests that numerical schemes with a small stencil, such as compact finite difference schemes, are preferred for the correction function method.

To locally solve the Cauchy problem \eqref{eqn:cau-prb}, we approximate the local solution in the narrow band $\Omega_{\Gamma}$ with a partition of unity approach.
Let the quasi-uniform point set $\left\{\boldsymbol{p}_i\right\}_{i=1}^{N_p} \subset\Gamma$ consist of primary points on the boundary $\Gamma$.
Let $\Omega_{\Gamma,i}$ be a neighborhood of the point $\boldsymbol{p}_i$.
Define $\Omega_{\Gamma}$ as the union of the neighborhoods
\begin{equation}
    \Omega_{\Gamma} = \bigcup_{i = 1}^{N_{p}} \Omega_{\Gamma,i}.
\end{equation}
Then $\{\Omega_{\Gamma, i}\}_{i=1}^{N_p}$ forms an overlapping decomposition of $\Omega_{\Gamma}$.
Note that each $\Omega_{\Gamma,i}$ should be chosen such that $\Omega_{\Gamma}$ covers all irregular grid nodes.
Unlike the original CFM \cite{Marques2011}, where $\Omega_{\Gamma}$ is defined as some particular grid patches relying on the cut pattern of $\Gamma$ with grid cells, the current definition of $\Omega_{\Gamma}$ is flexible since it only depends on the location of surface points.
This decomposition gives us a simple way to represent $C(\boldsymbol{x})$ in $\Omega_{\Gamma}$.

For the partitions $\Omega_{\Gamma, i}, i = 1, 2, \cdots, N_p$, define the compactly supported weight functions $\omega_i(\boldsymbol{x})$ such that $\text{supp}(\omega_i) = \Omega_{\Gamma,i}$ and 
\begin{equation}
    \sum_{i=1}^{N_p} \omega_i(\boldsymbol{y}) \equiv 1, \quad  \boldsymbol{y} \in \Omega_{\Gamma}=\bigcup_{i = 1}^{N_{p}} \Omega_{\Gamma,i}.
\end{equation}
In practice, the weight function $\omega_i$ can be constructed in many ways, such as Shepard’s method \cite{SHEPARDD1968}.
In this work, we use a simple non-smooth weight function,
\begin{equation}
    \omega_i(\boldsymbol{x}) = 
    \left\{
    \begin{aligned}
    1, &\quad \text{if $\boldsymbol{p}_i$ is the closest point to $\boldsymbol{x}$ for $i=1, 2, \cdots, N_p$,}\\
    0,& \quad \text{otherwise}.
    \end{aligned}
    \right.
\end{equation}
We remark that the smoothness of the weight function has a negligible effect on the algorithm.
The above simple weight function works very well for all numerical experiments.
Suppose $C_{h,i}(\boldsymbol{x})$ is an approximation to $C(\boldsymbol{x})$ for $\boldsymbol{x}\in\Omega_{\Gamma,i}$.
With the partition of unity, the complete approximate solution $C_h(\boldsymbol{x})$ for $\boldsymbol{x}\in\Omega_{\Gamma}$ is constructed as a linear combination of local solutions $C_{h,i}$,
\begin{equation}
    C_h(\boldsymbol{x}) = \sum_{i=1}^{N_p} \omega_i (\boldsymbol{x}) C_{h,i} (\boldsymbol{x}), \quad \boldsymbol{x}\in\Omega_{\Gamma}.
\end{equation}

To this end, we restrict the Cauchy problem \eqref{eqn:cau-prb} to the partition $\Omega_{\Gamma, i}$ and consider numerically solving a sequence of subproblems for $i=1,2,\cdots N_p$,
\begin{equation}\label{eqn:local-cau-prb}
    \begin{aligned}
    \Delta C_i(\boldsymbol{x}) - \kappa C_i(\boldsymbol{x})&= \Tilde{f}(\boldsymbol{x}),\quad &\boldsymbol{x}\in\Omega_{\Gamma,i},\\
    C_i(\boldsymbol{x}) &= a(\boldsymbol{x}),\quad &\boldsymbol{x}\in\Gamma\cap \Omega_{\Gamma,i},\\
    \partial_{\boldsymbol{n}} C_i(\boldsymbol{x}) &= b(\boldsymbol{x}), \quad &\boldsymbol{x}\in \Gamma\cap\Omega_{\Gamma,i}.
    \end{aligned}
\end{equation}
to obtain numerical solutions $C_{h,i}(\boldsymbol{x})$.
The restricted problems \eqref{eqn:local-cau-prb} are both temporally and spatially local, which explains the terminology ``local Cauchy problem''.

The method is more understandable if one regards the normal direction of $\Gamma$ as a time variable and the problems \eqref{eqn:local-cau-prb} as initial-boundary value problems (IBVPs).
Solving the restricted problems for the full Cauchy problem \eqref{eqn:cau-prb} resembles the explicit method for time-dependent PDEs.
In the correction function method, one does not need to be concerned with the stability of the explicit method since the solution is computed only one step away from the boundary $\Gamma$.
\begin{rem}
\label{rem:f-smooth}
    To obtain an accurate correction function $C(\boldsymbol{x})$, the right-hand side $\Tilde{f}(\boldsymbol{x})$ should be sufficiently smooth.
    For a fourth-order method, $C(\boldsymbol{x})$ is required to be at least $C^4$, and, consequently, $\Tilde{f}(\boldsymbol{x})$ is required to be at least $C^2$.
    Numerically, we can use the same partition of unity approach to represent $\Tilde{f}$ in $\Omega_{\Gamma}$.
    In each $\Omega_{\Gamma,i}$, $\Tilde{f}$ is replaced by a simple quadratic function using the jump information of $f$ (for example, in 2D, we use $[f]$, $[f_x]$, $[f_y]$, $[f_{xx}]$, $[f_{yy}]$, and $[f_{xy}]$).
    There are also several different ways to obtain smooth $f^+$ and $f^-$, such as the PDE-based method \cite{Aslam2004} and the partition of unity extension (PUX) method \cite{Fryklund2018}.
\end{rem}
\subsubsection{A mesh-free collocation method}

Let $\{\phi_{l,m,n}(\boldsymbol{x})\}_{l+m+n \leq p}$ denote the basis of Taylor polynomials of degree no more than $p$, where the subscripts $l$, $m$ and $n$ are non-negative integers.
The elements of the basis are given by, for example,
\begin{equation}
    \begin{aligned}
        &\phi_{0,0,0}(x,y,z) = 1,\\
        &\phi_{1,0,0}(x,y,z) = x,\quad \phi_{0,1,0}(x,y,z) = y,\quad \phi_{0,0,1}(x,y,z) = z,\\
        &\phi_{2,0,0}(x,y,z) = x^2,\quad \phi_{0,2,0}(x,y,z) = y^2,\quad \phi_{0,0,2}(x,y,z) = z^2,\\
        &\phi_{1,1,0}(x,y,z) = xy,\quad \phi_{1,0,1}(x,y,z) = xz,\quad \phi_{0,1,1}(x,y,z) = yz,\\
        &\cdots.
    \end{aligned}
\end{equation}
The approximate solution $C_{h,i}(\boldsymbol{x})$ is expressed as a linear combination of the basis functions
\begin{equation}
    C_{h,i}(\boldsymbol{x}) = \sum_{l+m+n \leq p} d_{l,m,n} \phi_{l,m,n}(\xi, \eta, \zeta) ,\quad \boldsymbol{x}\in\Omega_{\Gamma,i},
\end{equation}
where $\xi$, $\eta$ and $\zeta$ are scaled local coordinates of $\boldsymbol{x}$.
Suppose $\boldsymbol{p}_i = (x^{(i)}, y^{(i)}, z^{(i)})$ is the center point of the local domain $\Omega_{\Gamma,i}$.
The scaled local coordinate $\Tilde{\boldsymbol{x}} = (\xi,\eta,\zeta)$ of $\boldsymbol{x} = (x, y, z)$ is defined as
\begin{equation}
    \xi = (x - x^{(i)}) / h,\quad \eta = (y - y^{(i)}) / h,\quad \zeta = (z - z^{(i)}) / h,
\end{equation}
where $h$ is the mesh parameter.
To determine the coefficients $d_{l,m,n}$, we replace $C_i$ with $C_{h,i}$  in the problem \eqref{eqn:local-cau-prb} and let the equations be exactly satisfied at multiple points.
The resulting method is essentially mesh-free and falls into the category of collocation methods.
Then the chosen points are called ``collocation points.''
Since the problem \eqref{eqn:local-cau-prb} involves both the bulk PDE and boundary conditions, it involves collocation points in both $\Omega_{\Gamma,i}$ and $\Gamma\cap\Omega_{\Gamma,i}$.
Collocation points can be classified into three types based on the equations at which they are satisfied.
Let $\boldsymbol{x}_j^{pde}, j = 1, 2,\cdots, m_1$ be the points in $\Omega_{\Gamma,i}$ where the PDE is satisfied. 
Let $\boldsymbol{x}_j^{D}, j = 1, 2,\cdots, m_2$ and $\boldsymbol{x}_j^{N}, j = 1, 2,\cdots, m_3$ be the points on $\Gamma\cap\Omega_{\Gamma,i}$ where the Dirichlet and Neumann conditions are satisfied, respectively.
The problem \eqref{eqn:local-cau-prb} is approximated by the finite-dimensional problem
\begin{equation}\label{eqn:cau-appro}
    \begin{aligned}
    \sum_{l+m+n \leq p} (\Delta - \kappa)\phi_{l,m,n}(\Tilde{\boldsymbol{x}}^{pde}_j) d_{l,m,n} &= \Tilde{f}(\boldsymbol{x}_j^{pde}),\quad &\text{for } j = 1,2,\cdots,m_1 \\
    \sum_{l+m+n \leq p} \phi_{l,m,n}(\Tilde{\boldsymbol{x}}^{D}_j) d_{l,m,n} &= a(\boldsymbol{x}_j^{D}), \quad&\text{for } j = 1,2,\cdots,m_2, \\
    \sum_{l+m+n \leq p} \boldsymbol{n}(\boldsymbol{x}^{N}_j)\cdot \nabla\phi_{l,m,n}(\Tilde{\boldsymbol{x}}^{N}_j) d_{l,m,n} &= b(\boldsymbol{x}_j^{N}), \quad &\text{for } j = 1,2,\cdots,m_3.
    \end{aligned}
\end{equation}
The approximate problem \eqref{eqn:cau-appro} forms a linear system 
\begin{equation}\label{eqn:cau-ls}
    \mathbf{M}\mathbf{U} = \mathbf{Q},
\end{equation}
where the unknown vector $\mathbf{U}$ consists of the coefficients $d_{l,m,n}$.

\begin{rem}
The collocation method is closely related to the local  coordinate-transformation approach used in previous works \cite{Ying2007,Ying2013,Ying2014,Xie2020}.
The coordinate-transformation approach can also be viewed as a method for solving the local Cauchy problem \eqref{eqn:local-cau-prb} since the correction function $C(\boldsymbol{x})$ can also be approximated with the derivative jumps $[u], [u_x], [u_y], [u_z], [u_{xx}] \cdots$ in terms of a Taylor polynomial.
However, the derivation of derivative jumps in the coordinate-transformation approach involves repeatedly taking tangential derivatives and applying the chain rule, which requires tedious calculation, especially for high-order and 3D cases.
The collocation method introduced here is much simpler since applying the chain rule is not required.
\end{rem}

\subsection{Selection of collocation points}
Selecting collocation points is an essential part of the mesh-free collocation method to ensure accuracy and stability of the algorithm.
Different selection procedures for collocation points result in different systems \eqref{eqn:cau-ls} and different results.
For example, one can choose many collocation points such that their number is much more than the number of unknowns.
In that case, the linear system \eqref{eqn:cau-ls} becomes overdetermined and can be solved in the least-square sense, which is similar to the method in \cite{Marques2011}.
Here, an interpolation-type method is employed so that each equation in the system \eqref{eqn:cau-ls} is accurately satisfied.
An advantage of using an interpolation-type method is that when a boundary point coincides with a grid node, the correction function is accurate at the point since the Dirichlet jump condition $[u]$ is enforced accurately.

Before describing the selection procedure of collocation points, we emphasize a few key rules:
\begin{enumerate}
    \item[(a)] Collocation points should be chosen in $\Omega_{\Gamma,i}$ for the PDE and on $\Gamma\cap \Omega_{\Gamma, i}$ for boundary conditions.
    \item[(b)] Collocation points of the same type should be well-separated such that the resulting linear system is non-singular. 
    \item[(c)]  For each equation in \eqref{eqn:local-cau-prb}, the number of collocation points should be chosen to meet the formal accuracy requirement.
\end{enumerate}

Rule (a) is a basic requirement for consistency of the collocation method.
Rule (b) is intended to avoid a nearly singular or rank-deficient matrix $\mathbf{M}$ and to ensure the stability of the method.
For collocation points of the same type to be well-separated, the distance between two different points should have a positive lower bound.
Moreover, the number of projections of these points onto each spatial direction should be sufficiently large such that the interpolation bases associated with the points can span the polynomial space.
Rule (c) ensures accuracy of the collocation method. 
Note that the three equations in \eqref{eqn:local-cau-prb} have different orders of derivatives of $C_i$, and thus a polynomial approximation of $C_i$ results in different orders of accuracy for each equation. 
Since the equations in \eqref{eqn:local-cau-prb} are enforced accurately at collocation points, the collocation problem is also referred to as an interpolation problem.
With the error estimation of polynomial interpolation, one can find that the approximation errors at a point $\boldsymbol{x}$, away from collocation points, satisfy
\begin{equation}\label{eqn:appro-prob}
    \begin{aligned}
    (\Delta - \kappa) C_{h,i}(\boldsymbol{x}) - \Tilde{f}(\boldsymbol{x}) &= \mathcal{O}(h^{p-1}) ,\quad &\boldsymbol{x}\in\Omega_{\Gamma,i},\\
    C_{h,i}(\boldsymbol{x}) - a(\boldsymbol{x})&= \mathcal{O}(h^{p+1}),\quad &\boldsymbol{x}\in\Gamma\cap \Omega_{\Gamma,i},\\
    \partial_{\boldsymbol{n}} C_{h,i}(\boldsymbol{x}) - b(\boldsymbol{x})&= \mathcal{O}(h^{p}), \quad &\boldsymbol{x}\in \Gamma\cap\Omega_{\Gamma,i}.
    \end{aligned}
\end{equation}

\begin{figure}[htbp]
    \centering
    \subfigure[]{\includegraphics[width=0.25\textwidth]{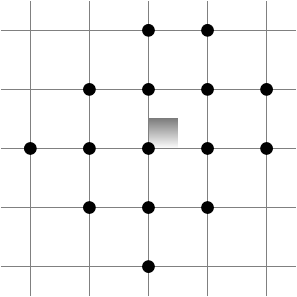}}\hspace{.2in}
    \subfigure[]{\includegraphics[width=0.25\textwidth]{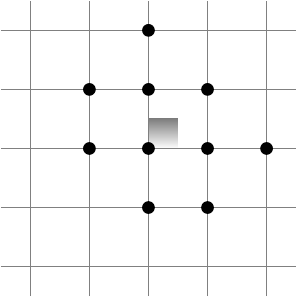}}\hspace{.2in}
    \subfigure[]{\includegraphics[width=0.25\textwidth]{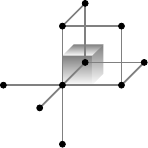}}
    \caption{Schematics of collocation points in 3D for (a) the Dirichlet boundary condition; (b) the Neumann boundary condition; and (c) the PDE. For center points of the local Cauchy problem that are located in the shaded region, collocation points are marked as black circles. Figures (a) and (b) show the projections of collocation points on the reference plane.}
    \label{fig:interp-stc}
\end{figure}

To take into account the consistency and stability requirements and to balance the approximation errors, we choose collocation points as interpolation points such that the corresponding Lagrange interpolant on these points has the same order of accuracy as shown in \eqref{eqn:appro-prob}.
Precisely, collocation points are chosen as interpolation points of a polynomial of degree (i) $(p-2)$ for the PDE; (ii) $p$ for the Dirichlet boundary condition; and (iii) $(p-1)$ for the Neumann boundary condition.
It should be mentioned that the Lagrange interpolant associated with the PDE is in $d$ space dimensions and those for the boundary conditions are in $(d-1)$ space dimensions.
Therefore, to choose collocation points for boundary conditions, we first project the boundary $\Gamma$ into its reference plane locally such that we can find the local stencil by working with the Cartesian grid on the planar domain. 
A good choice of the distribution of collocation points is illustrated in \Cref{fig:interp-stc}.
Similar point selection strategies for multivariate interpolation are used in \cite{Xie2020,Ying2014,Ying2013}.

If collocation points are chosen as above, the number of collocation points equals the number of degrees of freedom.
For example, in three space dimensions, the numbers of collocation equations for the PDE and the Dirichlet and Neumann boundary conditions are $\sum_{i=0}^{p-2} (i+1)(i+2)/2$, $(p+1)(p+2)/2$ and $p(p+1)/2$, respectively.
Obviously, it yields
\begin{equation}
    N^{eqn} = \sum_{i=0}^{p-2}\dfrac{(i+1)(i+2)}{2}+\dfrac{(p+1)(p+2)}{2}+\dfrac{p(p+1)}{2} = \sum_{i=0}^{p} \dfrac{(i+1)(i+2)}{2} = N^{dof}.
\end{equation}
Then the system \eqref{eqn:cau-ls} is a square one.
One can easily verify that similar results hold for the two-dimensional case as well.
The invertibility of the matrix $\mathbf{M}$ is difficult to prove since it depends on the geometry of $\Gamma$.
Nevertheless, if the collocation points are chosen as aforementioned, the linear system is always uniquely solvable with a standard decomposition method, such as the QR decomposition method.

\begin{rem}
We suggest using the scaled local coordinate $(\xi,\eta,\zeta)$ instead of the original coordinate $(x,y,z)$ for solving the problem \eqref{eqn:cau-appro}.
It is equivalent to rescaling the local Cauchy problem such that its characteristic length changes from $\mathcal{O}(h)$ to $\mathcal{O}(1)$.
Thus, the condition number of the problem \eqref{eqn:cau-ls} is essentially independent of the grid size $h$. The scaling can improve the accuracy and stability of the algorithm by reducing the effect of round-off error.
In the numerical experiments, the condition number $\text{cond}(\mathbf{M})$ is always on the order of $10^2 \sim 10^3$ regardless of how small the grid size is.
\end{rem}

\subsection{Extracting boundary data}

After solving the linear system of the corrected finite difference scheme, one can obtain the numerical solution at Cartesian grid nodes.
However, in the KFBI method, one needs to frequently use boundary/interface data, such as boundary value or normal derivative of the solution, at boundary nodes rather than Cartesian grid nodes.
In order to extract boundary data of the numerical solution, Lagrange interpolation is used to compute off-grid data.
One should also take into account the jump values of the potential function such that the Lagrange interpolation has high-order accuracy.
The correction function $C(\boldsymbol{x})$ introduced before now offers a suitable way to take into account the non-smoothness of the solution.
With the correction function, it is simple to reconstruct smooth data for interpolation using the piecewise smooth grid value.

For example, given a boundary point $\boldsymbol{p}\in\Gamma$, we try to obtain the one-sided limit boundary data of the numerical solution $v_h$ in $\Omega^+$.
Let $\boldsymbol{q}_i, i = 1, 2, \cdots$ be the grid nodes in the interpolation stencil near $\boldsymbol{p}$.
Suppose the numerical solution $v_h$ is piecewise smooth and coincides with the smooth functions $v_h^+$ and $v_h^-$ in $\Omega^+$ and $\Omega^-$, respectively.
We add the correction function $C(\boldsymbol{q}_i)$ to the grid value $v_h(\boldsymbol{q}_i)$ if $\boldsymbol{q}_i\in\Omega^-$ so that the interpolation data are smooth.
A Taylor expansion at $\boldsymbol{p}$ yields
\begin{align}
    v_h(\boldsymbol{q}_i) &= \sum_{l+m+n\leq p}\dfrac{p!}{l!m!n!}  \xi^l\eta^m\zeta^n \dfrac{\partial^{l+m+n}}{\partial x^l \partial y^m \partial z^n} v_h^+(\boldsymbol{p}) + \mathcal{O}(|\boldsymbol{q}_i-\boldsymbol{p}|^{p+1}),\quad  \text{if } \boldsymbol{q}_i \in \Omega^+,\\
    v_h(\boldsymbol{q}_i) + C(\boldsymbol{q}_i) &= \sum_{l+m+n\leq p}\dfrac{p!}{l!m!n!}  \xi^l\eta^m\zeta^n \dfrac{\partial^{l+m+n}}{\partial x^l \partial y^m \partial z^n} v_h^+(\boldsymbol{p}) + \mathcal{O}(|\boldsymbol{q}_i-\boldsymbol{p}|^{p+1}), \quad \text{if } \boldsymbol{q}_i \in \Omega^-.
\end{align}
where $(\xi, \eta, \zeta)^T = \boldsymbol{q}_i - \boldsymbol{p}$.
Now, by solving the interpolation problem, the function value and derivatives of $v_h^+(\boldsymbol{p})$ are obtained.


\section{Algorithm Summary} \label{sec:algo}
In this section, we summarize the proposed method.
We take the boundary integral equation \eqref{eqn:bie-ifp-1} as an example.
The algorithms for the boundary integral equations \eqref{eqn:bie-D},\eqref{eqn:bie-N},\eqref{eqn:bie-ifp-2} are similar.
The algorithm is summarized in \Cref{alg:IFP}.
\begin{algorithm}[htbp]
\caption{Correction function-based KFBI method}
\label{alg:IFP}
\begin{enumerate}[Step 1]
    \item Compute the right-hand side of \eqref{eqn:bie-ifp-1}, in which the integral operators $\mathcal{S}_i, \mathcal{K}_i, \mathcal{K}^{\prime}_i, \mathcal{D}_i, \mathcal{G}_i,\partial_{\mathbf{n}}\mathcal{G}f, i = 1,2$ are computed using the same approach in Step 3;
    \item Give an initial guess for $\varphi$ and $\psi$;
    \item Compute the integral operators $\mathcal{S}_i, \mathcal{K}_i, \mathcal{K}^{\prime}_i, \mathcal{D}_i, i = 1,2$ by solving the equivalent interface problem \eqref{eqn:equi-ifp};
    \begin{enumerate}[Step 3.1]
        \item Compute the correction function in $\Omega_{\Gamma}$ by solving the local Cauchy problem \eqref{eqn:cau-prb};
        \item Compute the correction terms in the right-hand side of \eqref{eqn:corrected-FDS};
        \item Solve the linear system of the finite difference scheme \eqref{eqn:corrected-FDS} with FFT;
        \item Compute the integral operator values by interpolation from the grid solution;
    \end{enumerate}
    \item Generate the next $\varphi$ and $\psi$ using the GMRES method and repeat Step 3 until the residual is less than a given tolerance.
\end{enumerate}

\end{algorithm}

In each iteration, individually computing the integral operators would require a total of eight calls of the FFT solver.
We stress that the number can be reduced to two since the terms $\mathcal{K}_i \varphi-\mathcal{S}_i \psi$ and $\mathcal D_i \varphi-\mathcal{K}^{\prime}_i \psi$ for $i=1$ or $2$ can be computed by calling the FFT solver only once.
By the principle of linear superposition, one only needs to solve the interface problem for the potential $D\varphi -S\psi$ and interpolate the function value and normal derivative on $\Gamma$ to obtain the terms.
In this way, only two calls of the FFT solver are required in each GMRES iteration.

\section{Numerical results}\label{sec:res}
In this section, numerical results for boundary and interface problems in both two and three space dimensions are presented.
In the following examples, irregular domains and interfaces are given in their level-set forms, which will be specified for each case.
Irregular domains and interfaces are embedded into a bounding box $\mathcal{B}$,  which is chosen as a square in 2D and a cube in 3D.
The box $\mathcal{B}$ is uniformly partitioned into $N$ intervals in each direction for simplicity.
The total number of primary boundary points representing the interface $\Gamma$ is denoted by $N_b$.

The following numerical experiments are performed on a personal computer with a 3.80 GHz Intel Core i7 processor. 
The codes for conducting the numerical experiments are written in C++. 
The tolerance in the GMRES method is fixed at $10^{-10}$.
GMRES iteration numbers and CPU times (in seconds) are reported to quantify the computational complexity.
Numerical errors on the grid node set $\Omega^h$ in the $L_2$ and maximum norms are defined as
\begin{equation}
    \Vert e\Vert_{2} = \sqrt{\dfrac{\sum_{\boldsymbol{x}\in\Omega^h}|v(\boldsymbol{x}) - u(\boldsymbol{x})|^2}{M}}, \quad \Vert e\Vert_{\infty} = \max_{\boldsymbol{x}\in\Omega^h}|v(\boldsymbol{x})-u(\boldsymbol{x})|,
\end{equation}
where $M$ is the number of grid nodes in $\Omega^h$, and $v$ and $u$ are the numerical and exact solutions, respectively.

\subsection{Two-dimensional examples}
\subsubsection{Boundary value problem}
In the first example, we solve the 2D Dirichlet BVP of the Poisson equation on a rotated ellipse-shaped domain $\Omega$
\begin{equation}
    \Omega = \left \{(x,y)\in\mathbb{R}^2: \dfrac{(x\cos\theta+y\sin\theta)^2}{a^2} + \dfrac{(y\cos\theta-x\sin\theta)^2}{b^2} < 1 \right\},
\end{equation}
with $a = 1, b = 0.5, \theta = -\pi/6$. 
The ellipse is embedded into the bounding box $\mathcal{B} = [-1.2,1.2]^2$.
The boundary condition and right-hand side are taken such that the exact solution satisfies
\begin{equation}
    u(x,y) = \exp(x)\sin(\cos(\pi/3) x + \sin(\pi/3) y).
\end{equation}

Numerical results are summarized in \Cref{tab:2d-poi-o4}. 
Nearly fifth-order accuracy in both the $L_2$ and maximum norms can be observed.
The increase in convergence order may be caused by the error of quartic polynomial interpolation, which is fifth-order accurate and dominates the numerical error in the vicinity of the boundary.
As the grid refines, the GMRES iteration number is essentially independent of grid size, which is a main advantage of the present method.
Taking into account FFT solvers and boundary operations in each iteration, the overall computational complexity of the method is given by $\mathcal{O}(N^2\log N + N_b)$ in two space dimensions.
On coarse grids, the CPU time scaling is close to $\mathcal{O}(N_b)$, implying that boundary operations dominate the computational cost.
On finer grids, the CPU time is roughly linearly proportional to $\mathcal{O}(N^2\log N)$, which implies that the computational cost is dominated by the FFT solver.
Isocontours of the numerical solution are also presented in \Cref{fig:poi-ell}.

\begin{table}[htbp]
\centering
\caption{Numerical results for the Dirichlet BVP of the Poisson equation on an ellipse-shaped domain.}
\label{tab:2d-poi-o4}
\begin{tabular}{|c|c|c|c|c|c|}
\hline
grid size                & 64$\times$64 & 128$\times$128 & 256$\times$256 & 512$\times$512 & 1024$\times$1024 \\ \hline
$N_{b}$                  & 116          & 230            & 460            & 918            & 1838             \\ \hline
itr no.                  & 10           & 10             & 9              & 9              & 9                \\ \hline
$\Vert e \Vert_{2}$      & 7.40E-06     & 1.12E-07       & 3.03E-09       & 6.86E-11       & 2.31E-12         \\ \hline
$\Vert e \Vert_{\infty}$ & 1.31E-04     & 3.69E-06       & 1.03E-07       & 3.56E-09       & 1.24E-10         \\ \hline
CPU time                 & 3.91E-03     & 6.35E-03       & 1.86E-02       & 5.81E-02       & 2.33E-01         \\ \hline
\end{tabular}
\end{table}

\subsubsection{Interface problem with multiple interfaces}
In the second example, we solve the 2D Poisson interface problem with multiple disjoint interfaces, which are eight circles and a five-fold star, on the domain $\mathcal{B} = [-1.7,1.7]^2$.
The circles are given by 
\begin{equation}
    \Gamma_{m}^{cir} = \left\{ (x,y)\in\mathbb{R}^2: (x-\cos(m\pi/4))^2 + (y-\sin(m\pi/4))^2 = r^2 \right\},\quad m = 1, 2, \cdots, 8,
\end{equation}
with $r = 0.383$.
The five-fold star is given by
\begin{equation}
    \Gamma^{star} =\left\{(x,y)\in\mathbb{R}^2: \dfrac{x^2}{a^2} + \dfrac{y^2}{b^2} = (1.0 + \varepsilon \sin(m \arctan(\dfrac{y}{x})))^2\right\},
\end{equation}
with $a = b = 0.514, \varepsilon=0.2, m = 5$.
Two adjacent interfaces may become very close to each other, and, as a result, there may be more than one intersection point between two adjacent grid nodes.
The boundary condition, interface condition and right-hand side are chosen such that the exact solution is given by
\begin{equation}
    u(x, y) = \left\{
    \begin{aligned}
    &\exp(0.6x+0.8y), \quad & \text{in }\Omega_i,\\
    &\sin(\pi(x+1)/2)\sin(\pi(y+1)/2),\quad & \text{in }\Omega_e,
    \end{aligned}
    \right.
\end{equation}
where $\Omega_i$ denotes the union of the interiors of the circles and the star and $\Omega_e$ denotes the exterior domain.
The diffusion coefficients are chosen as $\sigma_i = 1$ in $\Omega_i$ and $\sigma_e = 3$ in $\Omega_e$.
For this and the following examples, the subscripts $i$ and $e$ represent variables in the interior and exterior regions, respectively.

\begin{figure}
\centering
  \begin{minipage}[t]{0.3\linewidth}
    \centering
    \includegraphics[width=\textwidth]{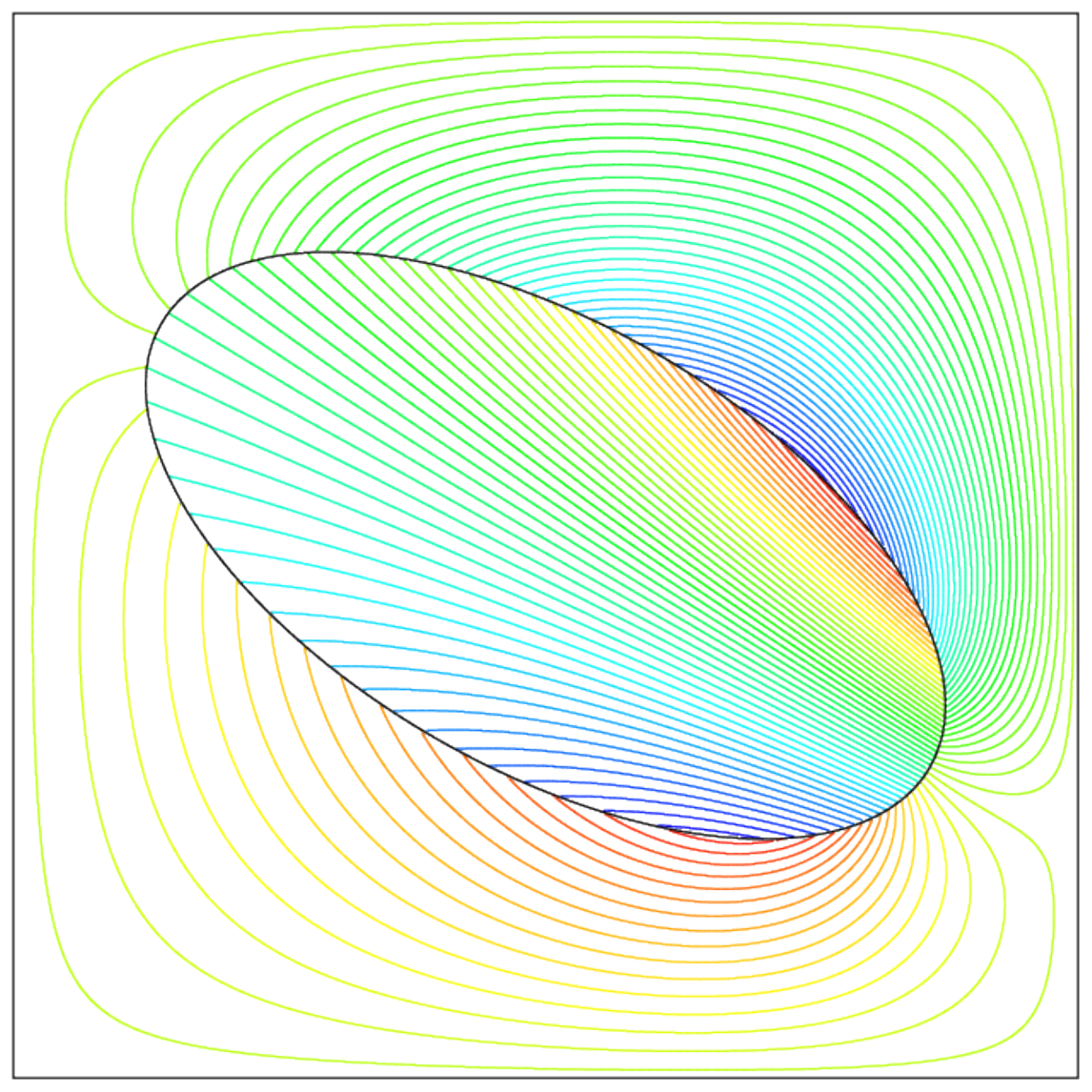}
    \caption{Numerical solution to the Dirichlet BVP of the Poisson equation on an ellipse-shaped domain.}
    \label{fig:poi-ell}
  \end{minipage}%
  \hspace{0.05\linewidth}
  \begin{minipage}[t]{0.3\linewidth}
    \centering
    \includegraphics[width=\textwidth]{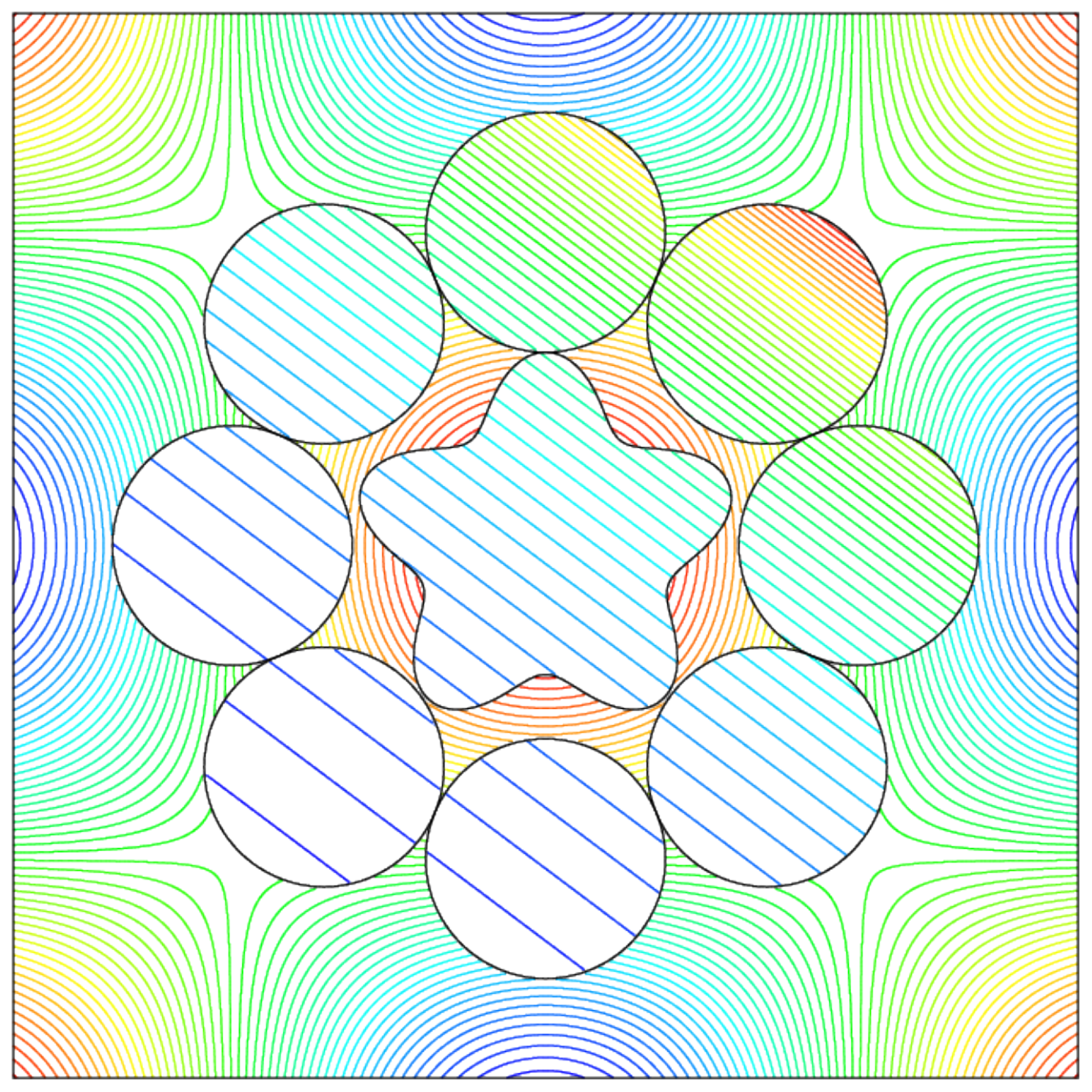}
    \caption{Numerical solution to the Poisson interface problem with multiple touching interfaces.}
    \label{fig:2d-poi-ifp-mul}
  \end{minipage}
\end{figure}

Numerical results are summarized in \Cref{tab:2d-poi-ifp-mul}.
The solutions in both the interior and exterior domains have fourth-order accuracy.
The GMRES iteration number is essentially independent of grid size, even if there are arbitrarily close interfaces.
It can be observed that the iteration number is slightly larger on the coarsest grid $N=64$.
A coarse Cartesian grid may not be able to accurately capture the geometry of complex interfaces.
This affects the well-conditioned property of the discrete boundary integral equation and causes the increase in iteration number.
Isocontours of the numerical solution are shown in \Cref{fig:2d-poi-ifp-mul}.

\begin{table}[htbp]
\centering
\caption{Numerical results for the Poisson interface problem with multiple touching interfaces.}
\label{tab:2d-poi-ifp-mul}
\begin{tabular}{|c|c|c|c|c|c|}
\hline
grid size                      & 64$\times$64 & 128$\times$128 & 256$\times$256 & 512$\times$512 & 1024$\times$1024 \\ \hline
$N_{b}$                        & 392          & 784            & 1568           & 3136           & 6282             \\ \hline
itr no.                        & 29           & 20             & 20             & 19             & 19               \\ \hline
$\Vert e \Vert_{\infty,\Omega_i}$ & 5.21E-05     & 4.65E-07       & 3.14E-09       & 1.90E-10       & 3.51E-11         \\ \hline
$\Vert e \Vert_{\infty,\Omega_e}$ & 3.94E-05     & 5.73E-07       & 3.19E-09       & 1.89E-10       & 3.51E-11         \\ \hline
CPU time                       & 7.81E-02     & 8.59E-02       & 1.48E-01       & 2.58E-01       & 5.86E-01         \\ \hline
\end{tabular}
\end{table}

\subsection{Three-dimensional examples}
To demonstrate the applicability of the present method, we consider solving three-dimensional problems.
\subsubsection{Poisson BVP}

This example is the Neumann BVP of the Poisson equation on a torus in 3D.
The torus is given by
\begin{equation}\label{eqn:ls-torus}
    \Omega = \left\{(x,y,z)\in\mathbb{R}^3: (1 - \sqrt{x^2+y^2})^2 + z^2 < 0.4^2\right\}.
\end{equation}
The bounding box is taken as $\mathcal{B}=[-1.5,1.5]^3$.
The boundary condition and right-hand side are taken such that the exact solution satisfies
\begin{equation}
    u(x,y,z) = \exp(z)(\cos(2x)+\cos(3y)).
\end{equation}
Note that the solution to the Poisson Neumann BVP is only determined up to an additive constant.
We first subtract a constant from the right-hand side of the linear system such that it has zero mean. 
At the same time, the matrix-vector products in the GMRES iterations are subtracted by a constant such that their means are zero.
To compute numerical errors, we need to add a constant to the numerical solution such that it matches the exact solution at a point.

Numerical results and the numerical solution are presented in \Cref{tab:3d-neu-poi} and \Cref{fig:poi-neu-3d}, respectively.
Fourth-order accuracy in both the $L_2$ and maximum norms is reached for the Neumann BVP. 
In this example, the GMRES iteration number decreases slightly as the grid refines.
Since the discrete linear system mimics the original well-conditioned BIE, the approximation with a fine grid is more accurate.
We believe that the better approximation property of a fine grid gives a linear system with a better condition number and is responsible for the faster convergence of the GMRES method.

Theoretically, the computational complexity in three space dimensions is $\mathcal{O}(N^3\log N + N_b)$.
The cost of boundary operations is more important than that in two space dimensions since the polynomial approximation for the correction function needs more terms in this case.
As a result, the overall computational cost is closer to $\mathcal{O}(N^2)$ since we have $N_b=\mathcal{O}(N^2)$.

\begin{table}[htbp]
\centering
\caption{Numerical results for the Neumann BVP of the Poisson equation on a torus.}
\label{tab:3d-neu-poi}
\begin{tabular}{|c|c|c|c|c|}
\hline
grid size & 64$\times$64$\times$64 & 128$\times$128$\times$128 & 256$\times$256$\times$256 & 512$\times$512$\times$512 \\ \hline
$N_b$                    & 6168     & 24656    & 98668    & 394548   \\ \hline
itr no.                  & 23       & 21       & 18       & 17       \\ \hline
$\Vert e \Vert_{2}$      & 2.39E-04 & 2.76E-05 & 2.38E-06 & 1.68E-07 \\ \hline
$\Vert e \Vert_{\infty}$ & 1.18E-03 & 7.97E-05 & 5.98E-06 & 4.01E-07 \\ \hline
CPU time                 & 1.62E+00 & 6.35E+00 & 2.50E+01 & 1.35E+02 \\ \hline
\end{tabular}
\end{table}

\subsubsection{Modified Helmholtz BVP}

\begin{figure}
\centering
  \begin{minipage}[t]{0.3\linewidth}
    \centering
    \includegraphics[width=\textwidth]{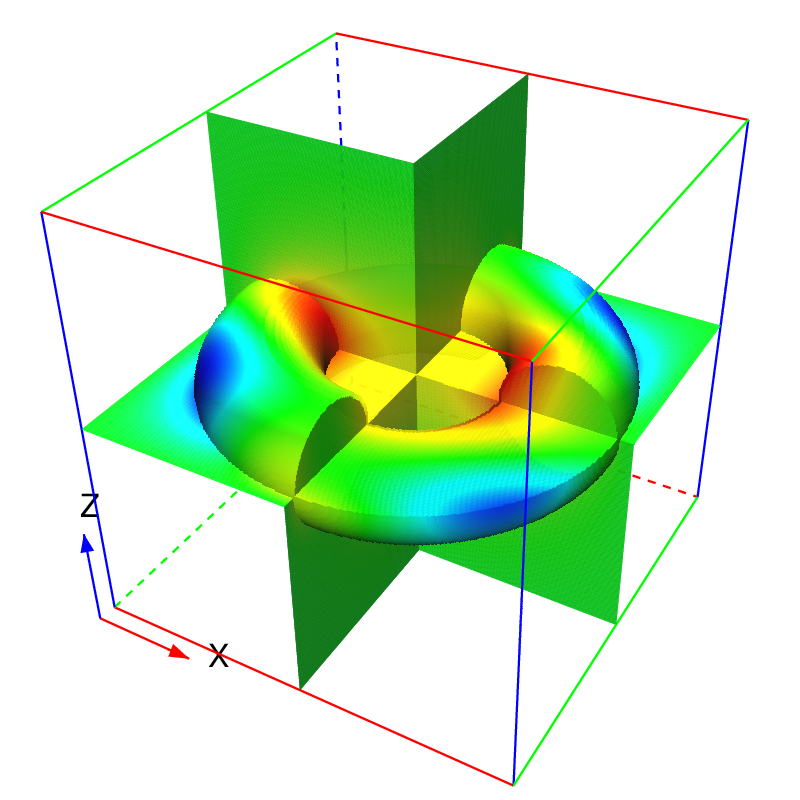}
    \caption{Numerical solution to the Neumann BVP of the Poisson equation on a torus.}
    \label{fig:poi-neu-3d}
  \end{minipage}%
  \hspace{0.05\linewidth}
  \begin{minipage}[t]{0.3\linewidth}
    \centering
    \includegraphics[width=\textwidth]{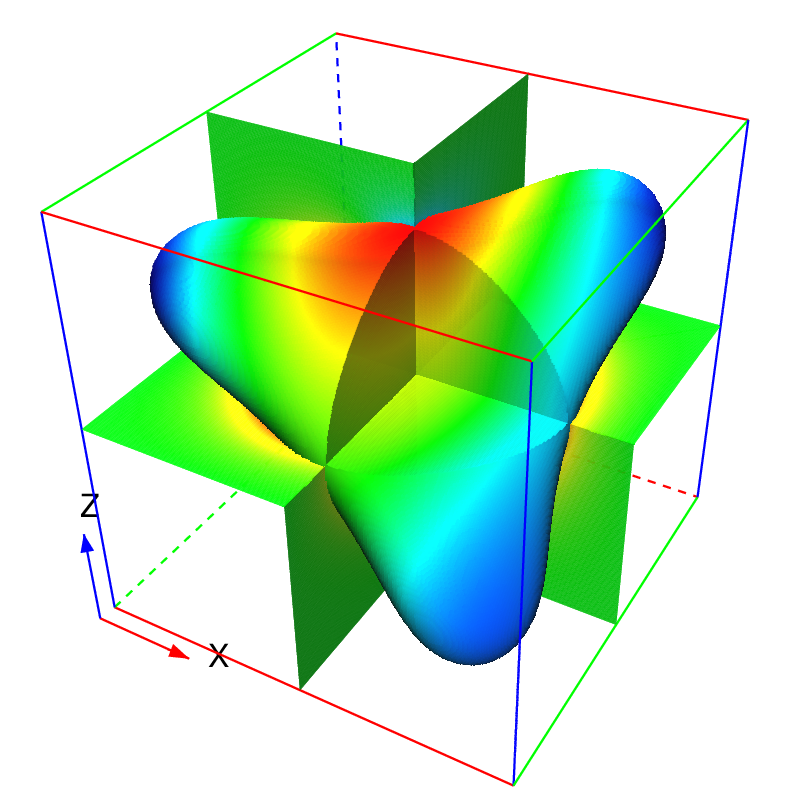}
    \caption{Numerical solution to the Dirichlet BVP of the modified Helmholtz equation.}
    \label{fig:hel-dir-3d}
  \end{minipage}
\end{figure}

As in the preceding example, we solve the Dirichlet BVP of the modified Helmholtz equation with $\kappa = 100$ on the domain $\Omega$, which is given by
\begin{equation}
    \Omega = \left\{(x,y,z)\in\mathbb{R}^3: (1+4x^2)(1+4y^2)(1+4z^2)+64xyz+4x^2+4y^2+4z^2<3 \right\}.
\end{equation}
This domain has relatively large curvature and is difficult to capture with a coarse grid.
The bounding box is taken as $[-0.7,0.7]^3$.
The boundary condition and right-hand side are chosen such that the exact solution satisfies
\begin{equation}
    u(x,y,z) = \exp(z)(\cos(5x)+\cos(2y)).
\end{equation}
Numerical results are summarized in \Cref{tab:3d-Diri-Hel}. 
The numerical solution is presented in \Cref{fig:hel-dir-3d}.
One can observe that the numerical error is large on the grid $N=64$ and decreases rapidly when the grid is refined to $N=128$.
It can be explained by the fact that the coarse grid $N=64$ may not be able to fully capture the fast changes of the boundary and cause large errors for near-interface corrections and surface interpolations.
As the grid refines, the decrease in numerical errors matches the fourth-order accuracy, as anticipated.
The coarse grid with $N=64$ also requires more GMRES iterations to converge. In each iteration, the CPU time scaling is close to $\mathcal{O}(N_b)$ due to the dominance of boundary operations.

\begin{table}[htbp]
\centering
\caption{Numerical results for the Dirichlet BVP of the modified Helmholtz equation. }
\label{tab:3d-Diri-Hel}
\begin{tabular}{|c|c|c|c|c|}
\hline
grid size & 64$\times$64$\times$64 & 128$\times$128$\times$128 & 256$\times$256$\times$256 & 512$\times$512$\times$512 \\ \hline
$N_b$                    & 7082     & 28482    & 114018   & 455450   \\ \hline
itr no.                  & 14       & 9        & 9        & 9        \\ \hline
$\Vert e \Vert_{2}$      & 1.14E-05 & 7.34E-08 & 2.82E-09 & 1.20E-10 \\ \hline
$\Vert e \Vert_{\infty}$ & 1.80E-03 & 6.34E-06 & 5.06E-07 & 2.17E-08 \\ \hline
CPU time                 & 2.98E+00 & 8.52E+00 & 3.82E+01 & 1.75E+02 \\ \hline
\end{tabular}
\end{table}

\subsubsection{Interface problem with high-contrast coefficients}
In this example, we solve the Poisson interface equation with a four-atom molecular-shaped interface in the domain $\mathcal{B} = [-1.2,1.2]^3$.
The interface $\Gamma$ is given by 
\begin{equation}
    \Gamma = \left\{\boldsymbol{x} = (x,y,z)\in\mathbb{R}^3: \sum_{k=1}^4 \exp(-\dfrac{|\boldsymbol{x} - \boldsymbol{x}_k|^2}{r^2}) = 0.6  \right\},
\end{equation}
with $\boldsymbol{x}_1 = (\sqrt{3}/3, 0,-\sqrt{6}/12)$, $\boldsymbol{x}_2 = (-\sqrt{3}/6, 0.5,-\sqrt{6}/12)$, $\boldsymbol{x}_3 = (-\sqrt{3}/6, -0.5,-\sqrt{6}/12)$ and $\boldsymbol{x}_4 = (0, 0,\sqrt{6}/4)$.

\begin{equation}
    u(x, y,z) = \left\{
    \begin{aligned}
    &\sin^2(2x)\cos^2(2y)\cos(z), \quad & \text{in the interior }\Omega_{i},\\
    &\cos(x)\cos(y)\cos(z),\quad & \text{in the exterior }\Omega_e.
    \end{aligned}
    \right.
\end{equation}
The coefficient ratio $\sigma_e/\sigma_i$ varies from $10$ to $10^4$, and its effect on the performance of the present method is studied in this example.
This effect was also studied by \cite{ZHOU20061,Wang2004,Marques2017}.
The numerical solution is shown in \Cref{fig:3d-ifp-ratio}.
According to the numerical results presented in \Cref{tab:3d-ifp-ratio}, high-contrast coefficients only have a small effect on the numerical accuracy, even for the extreme case $\sigma_e/\sigma_i=10^4$.
The GMRES iteration number is slightly affected by the coefficient ratio on coarse grids.
As the grid refines, the GMRES iteration number is rather stable and is independent of the coefficient ratio.
This is also due to the fact that a fine grid has a better approximation property, as aforementioned.

\begin{table}[htbp]
\centering
\caption{Numerical results for the Poisson interface problems with varying coefficient ratios. }
\label{tab:3d-ifp-ratio}
\begin{tabular}{|c|c|c|c|c|c|c|}
\hline
$\sigma_e:\sigma_i$     & N   & itr no. & $\Vert e \Vert_{2,\Omega_i}$ & $\Vert e \Vert_{\infty,\Omega_i}$ & $\Vert e \Vert_{2,\Omega_e}$ & $\Vert e \Vert_{\infty,\Omega_e}$ \\ \hline
\multirow{3}{*}{$10:1$} & 128 & 11      & 5.15E-08                     & 2.18E-07                          & 1.43E-08                     & 2.33E-07                          \\ \cline{2-7} 
                        & 256 & 10      & 3.25E-09                     & 1.22E-08                          & 9.03E-10                     & 1.33E-08                          \\ \cline{2-7} 
                        & 512 & 10      & 2.04E-10                     & 7.59E-10                          & 5.76E-11                     & 8.08E-10                          \\ \hline
\multirow{3}{*}{$10^2:1$} & 128 & 13      & 5.78E-08                     & 2.74E-07                          & 1.89E-08                     & 2.92E-07                          \\ \cline{2-7} 
                        & 256 & 11      & 3.63E-09                     & 1.61E-08                          & 1.18E-09                     & 1.71E-08                          \\ \cline{2-7} 
                        & 512 & 10      & 2.29E-10                     & 1.02E-09                          & 7.59E-11                     & 1.06E-09                          \\ \hline
\multirow{3}{*}{$10^4:1$} & 128 & 14      & 5.86E-08                     & 2.81E-07                          & 1.94E-08                     & 2.99E-07                          \\ \cline{2-7} 
                        & 256 & 11      & 3.68E-09                     & 1.66E-08                          & 1.21E-09                     & 1.75E-08                          \\ \cline{2-7} 
                        & 512 & 10      & 2.33E-10                     & 1.05E-09                          & 7.83E-11                     & 1.09E-09                          \\ \hline
\end{tabular}
\end{table}

\subsubsection{Interface problem with arbitrarily close interfaces}

\begin{figure}
\centering
  \begin{minipage}[t]{0.3\linewidth}
    \centering
    \includegraphics[width=\textwidth]{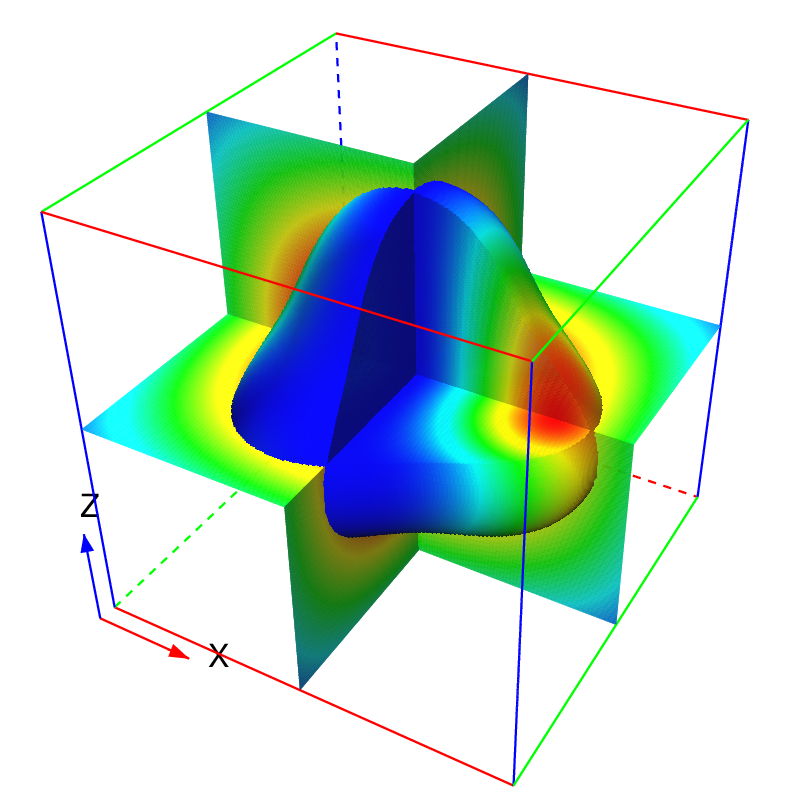}
    \caption{Numerical solution to the Poisson interface problem with varying coefficient ratios.}
    \label{fig:3d-ifp-ratio}
  \end{minipage}%
  \hspace{0.05\linewidth}
  \begin{minipage}[t]{0.3\linewidth}
    \centering
    \centering
    \includegraphics[width=\textwidth]{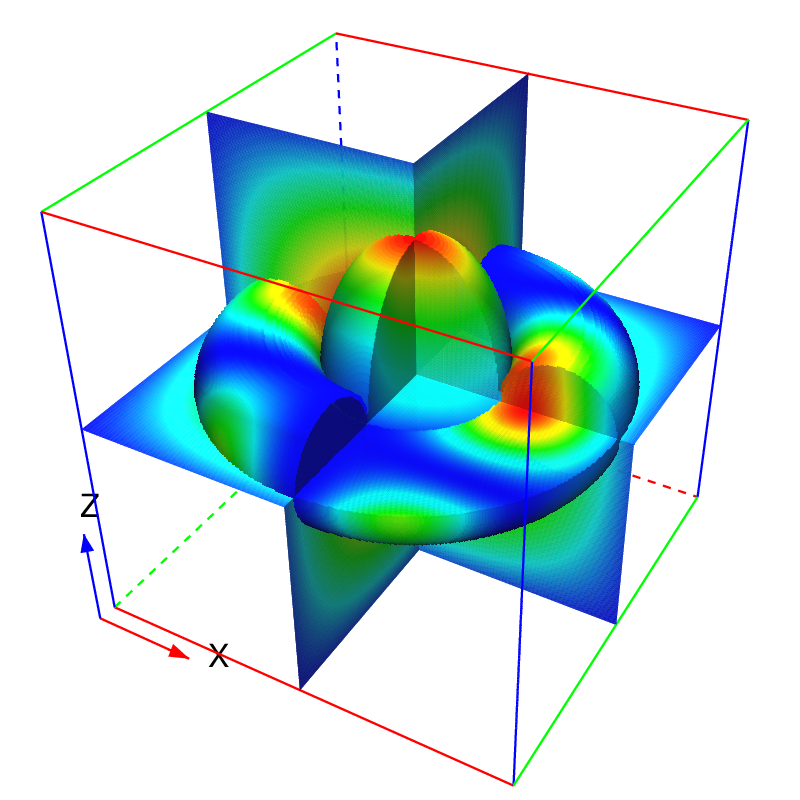}
    \caption{Numerical solution to the Poisson interface problem with a touching curve.}
    \label{fig:3d-poi-ifp}
  \end{minipage}
\end{figure}

In this case, we solve the Poisson interface problem with the presence of arbitrarily close interfaces in three space dimensions. 
Interfaces are taken as a torus and an ellipsoid. 
The torus-shaped interface $\Gamma^{tor}$ is given by the boundary of the domain $\Omega$ defined in \eqref{eqn:ls-torus}.
The ellipsoid-shaped interface is given by
\begin{equation}
    \Gamma^{ell} = \left\{(x,y,z)\in\mathbb{R}^3: \dfrac{x^2}{a^2} + \dfrac{y^2}{b^2}+\dfrac{z^2}{c^2} = 1\right\},
\end{equation}
with $a = b = 0.6, c = 1$.
The two interfaces are very close to each other near the curve
\begin{equation}
    S = \left\{(x,y,z)\in\mathbb{R}^3: \dfrac{x^2}{a^2} + \dfrac{y^2}{b^2} = 1, \quad z = 0\right\}.
\end{equation}
In this configuration, since the curve $S$ is a one-dimensional object, the number of multi-intersection grid line segments---grid line segments that intersect interfaces multiple times---is on the order of $\mathcal{O}(N)$.
The problem is challenging for classical body-fitted approaches because it is nearly impossible to resolve with a body-fitted mesh when the two interfaces are too close.
The bounding box $\mathcal{B}$ is taken as $[-1.5,1.5]^3$.
The boundary condition, interface condition and right-hand side are chosen such that the exact solution reads
\begin{equation}
    u(x, y,z) = \left\{
    \begin{aligned}
    &\sin^2(2x)\cos^2(2y)\cos(z), \quad & \text{in the torus }\Omega_{i, 1},\\
    &\exp(z)(\cos(2x)+\cos(3y)), \quad &\text{in the ellipsoid }\Omega_{i, 2},\\
    &\cos(x)\cos(y)\cos(z),\quad & \text{in the exterior region }\Omega_e.
    \end{aligned}
    \right .
\end{equation}
The coefficients are chosen as $\sigma_i = 1$ in $\Omega_{i, 1}\cup\Omega_{i, 2}$ and $\sigma_e = 3$ in $\Omega_e$.
Numerical results are summarized in \Cref{tab:poi-ifp-3d-toch-cur}. 
The numerical solution is visualized and shown in \Cref{fig:3d-poi-ifp}.
It is observed that fourth-order accuracy is achieved in all regions, except for an accuracy loss on the coarsest grid $N=64$ due to similar reasons that were mentioned before.

\begin{table}[htbp]
\centering
\caption{Numerical results for the Poisson interface problem with a touching curve.}
\label{tab:poi-ifp-3d-toch-cur}
\begin{tabular}{|c|c|c|c|c|}
\hline
grid size & 64$\times$64$\times$64 & 128$\times$128$\times$128 & 256$\times$256$\times$256 & 512$\times$512$\times$512 \\ \hline
$N_b$                     & 8738     & 34866    & 139454   & 557794   \\ \hline
itr no.                            & 22       & 21       & 19       & 17       \\ \hline
$\Vert e \Vert_{\infty,\Omega_{i,1}}$ & 2.33E-03 & 6.75E-07 & 4.38E-08 & 2.13E-09 \\ \hline
$\Vert e \Vert_{\infty,\Omega_{i,2}}$ & 5.56E-05 & 7.72E-07 & 4.39E-08 & 2.11E-09 \\ \hline
$\Vert e \Vert_{\infty,\Omega_e}$   & 3.24E-03 & 5.97E-07 & 3.65E-08 & 2.06E-09 \\ \hline
CPU time                         & 5.55E+00 & 2.08E+01 & 7.07E+01 & 2.86E+02 \\ \hline
\end{tabular}
\end{table}

\subsubsection{Heterogeneous interface problem}

\begin{figure}[htbp]
    \centering
    \includegraphics[width=0.3\textwidth]{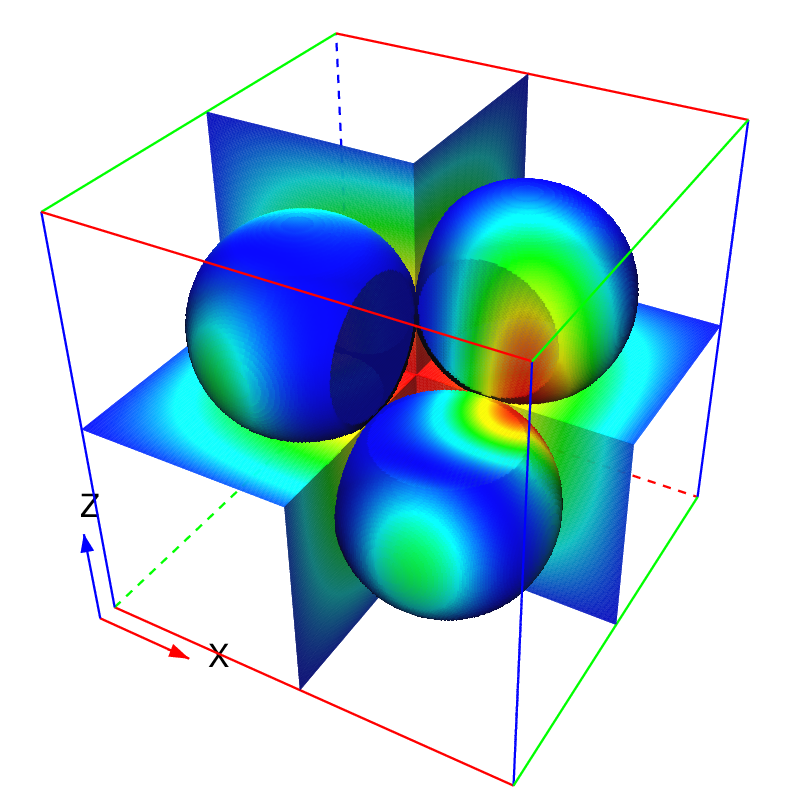}
    \caption{Numerical solution to the heterogeneous interface problem.}
    \label{fig:3d-ifp-heter}
\end{figure}

In the final example, we consider the heterogeneous interface problem in three space dimensions.
Interfaces are taken as three spheres with radius $r=0.7$ whose centers are chosen as $\boldsymbol{x}_1=(0.5,0.5,0.5)$, $\boldsymbol{x}_2=(-0.5,-0.5,0.5)$ and $\boldsymbol{x}_3=(0.5,-0.5,-0.5)$, respectively.
The coefficients on each side of the interfaces are given as 
\begin{equation}
    \sigma_i = 1,\quad \kappa_i = 0,\quad \sigma_e = 4,\quad \kappa_e = 10,
\end{equation}
such that the unknown function $u$ satisfies the Poisson equation in the interior region and the modified Helmholtz equation in the exterior region.
It is called a heterogeneous interface problem since the elliptic differential operators on the two sides of the interfaces are of different types.
The heterogeneous interface problem is a linearized version of the Poisson-Boltzmann equation, which appears in the Poisson-Boltzmann theory in biophysics for modeling solvated biomolecular systems.
\Cref{tab:3d-ifp-heter} and \Cref{fig:3d-ifp-heter} show the numerical results and the visualization of the numerical solution, respectively.
Once again we observe the fourth-order convergence in both regions.
The number of GMRES iterations is essentially independent of the grid size.

\begin{table}[htbp]
\centering
\caption{Numerical results for the heterogeneous interface problem.}
\label{tab:3d-ifp-heter}
\begin{tabular}{|c|c|c|c|c|}
\hline
grid size & 64$\times$64$\times$64 & 128$\times$128$\times$128 & 256$\times$256$\times$256 & 512$\times$512$\times$512 \\ \hline
$N_b$        & 6987     & 27939    & 111774   & 447264   \\ \hline
itr no.       & 21       & 20       & 20       & 20       \\ \hline
$\Vert e \Vert_{\infty,\Omega_i}$        & 6.71E-06 & 4.79E-07 & 3.32E-08 & 2.26E-09 \\ \hline
$\Vert e \Vert_{\infty,\Omega_e}$        & 3.89E-06 & 2.65E-07 & 1.73E-08 & 1.15E-09 \\ \hline
CPU time       & 1.34E+00 & 6.59E+00 & 4.08E+01 & 3.20E+02 \\ \hline
\end{tabular}
\end{table}

\section{Discussion}\label{sec:dis}
This work proposes a new version of the kernel-free boundary integral method for solving elliptic partial differential equations in two and three space dimensions with high accuracy.
The KFBI method solves boundary and interface problems with their boundary integral formulations.
It computes boundary and volume integrals by solving equivalent interface problems with fast PDE solvers and then obtains boundary values by interpolation.

The equivalent interface problems are simpler than the original problem and are essential for the KFBI method.
To accommodate the jump conditions across the interface, a correction function is introduced in the vicinity of the interface to derive corrected finite difference schemes and the boundary interpolation scheme.
Unlike the original KFBI method, which applies a local coordinate transformation to calculate correction terms, the new approach obtains correction terms by solving a local Cauchy problem for the correction function.
The local Cauchy problem is solved with a mesh-free collocation method, for which we also propose a strategy to choose collocation points such that the resulting linear system is accurate and stable.
The resulting method avoids repeatedly taking tangential derivatives of the jump conditions and significantly simplifies the derivation procedure.

The presented method is efficient and accurate, which is demonstrated through several challenging numerical experiments.
The efficiency of the method relies on the well-conditioning of the boundary integral equations and the applicability of fast PDE solvers (FFT and geometric multigrid methods) on a Cartesian grid.
Even though the presented numerical results are based on a fourth-order implementation of the method, the method can be extended to arbitrary accuracy in principle \cite{Marques2011}.

Finally, we emphasize that the present method is designed for implicitly defined interfaces with level-set formulations.
Although this work uses an analytic expression of the level-set function, extending the method to cases when the level-set function is only given at Cartesian grid nodes is straightforward.
It may have advantages for solving moving interface problems and free boundary problems when combined with the level-set method \cite{OSHER198812, Sethian1591}.

\section*{Acknowledgments}
This work is financially supported by the Shanghai Science and Technology Innovation Action Plan in Basic Research Area (Project No. 22JC1401700). It is also partially supported by the National Key R\&D Program of China (Project No. 2020YFA0712000), the Strategic Priority Research Program of the Chinese Academy of Sciences (Grant No. XDA25010405) and the National Natural Science Foundation of China (Grant No. DMS-11771290).

\bibliography{mybibfile}

\end{document}